\documentclass[10pt]{article}

\usepackage[letterpaper,top=0.5in, bottom=1in, left=0.8in, right=0.8in, headsep=0.in, centering]{geometry}
\usepackage[utf8]{inputenc}

\usepackage{amsmath}
\usepackage{bm}
\usepackage{amsfonts}
\usepackage{amssymb}
\usepackage{graphicx}
\usepackage{enumitem}
\usepackage{amsthm}
\usepackage{caption} 
\usepackage{array}
\usepackage{color}
\usepackage{csquotes}
\usepackage{bookmark}
\usepackage{float}
\usepackage{multirow}
\usepackage[normalem]{ulem} % delete line
\usepackage{orcidlink}

\MakeOuterQuote{"}

\makeatletter
\def\thickhline{%
  \noalign{\ifnum0=`}\fi\hrule \@height \thickarrayrulewidth \futurelet
   \reserved@a\@xthickhline}
\def\@xthickhline{\ifx\reserved@a\thickhline
               \vskip\doublerulesep
               \vskip-\thickarrayrulewidth
             \fi
      \ifnum0=`{\fi}}
\makeatother
\newlength{\thickarrayrulewidth}
\numberwithin{equation}{section}

\newcommand{\sd}{\mathord{\prime\mkern-2.5mu\reflectbox{$\scriptstyle\prime$}}}

\newcommand\innerprd[1][.5]{\mathbin{\vcenter{\hbox{\scalebox{#1}{$\bullet$}}}}}

\newtheorem{theorem}{Theorem}[section]
\newtheorem{lemma}[theorem]{Lemma}
\newtheorem{proposition}[theorem]{Proposition}
\newtheorem{corollary}[theorem]{Corollary}

\theoremstyle{definition}
\newtheorem{definition}[theorem]{Definition}
\newtheorem{example}[theorem]{Example}
\newtheorem{remark}[theorem]{Remark}

\title{
Regular specular differentiation in Euclidean spaces
}

\author{
Kiyuob Jung\,\orcidlink{0009-0007-5075-4061}\thanks{Department of Mathematics, Michigan State University, East Lansing, MI 48824, USA, E-mail: kyjung@msu.edu},
Jehan Oh\,\orcidlink{0000-0001-9602-5111}\thanks{Department of Mathematics, Kyungpook National University, Daegu, 41566, Republic of Korea, E-mail: jehan.oh@knu.ac.kr}}

\begin{document}
\date{}
\maketitle

\begin{abstract}
  We study the regular specular derivative, a generalized derivative defined at every point where both one-sided derivatives exist and are finite.
  Geometrically, it is the slope of the mirror that reflects the left tangent ray into the right one.
  In one variable we derive computational formulas, prove inverse function and rotation rules, establish Quasi-Rolle's Theorem and the Quasi-Mean Value Theorem, and obtain a derivative-limit theorem, which shows that twice regularly specularly differentiable functions are continuously differentiable.
  We also prove both parts of the Fundamental Theorem of Calculus.
  In several variables we introduce specular gradients, directional derivatives, tangent hyperplanes, and normal vectors, show that a continuous specular gradient forces classical differentiability, and characterize when the specular tangent hyperplane is unique.
\end{abstract}

\smallskip


\smallskip


\section{Introduction}

Derivatives measure how a function changes, and much of analysis is built on them.
The classical derivative, however, requires the two one-sided difference quotients to converge to a common limit.
This requirement fails at corners, which arise from absolute values, maxima, and minima, and hence in optimization and in many applied problems.
This limitation has motivated a variety of generalized derivatives, including symmetric derivatives, Dini derivatives, subderivatives, and weak derivatives; see \cite{1966_Bruckner,1994_Bruckner_BOOK} for extensive surveys.
Each generalization keeps some features of the classical derivative and gives up others.
Subderivatives keep the geometry of supporting lines but produce a set of slopes rather than a single number.
Weak derivatives keep the integration by parts formula but lose their pointwise meaning.
Throughout this paper, we call $f'$ the \emph{classical derivative} of $f$ to distinguish it from the generalized derivatives under discussion.

The closest relative of the object studied here is the symmetric derivative, which replaces the ordinary difference quotient by the quotient of $f(x+h)-f(x-h)$ over $2h$.
The symmetric derivative coincides with the classical derivative wherever the latter exists, but it may also exist at points where the classical derivative does not.
It need not satisfy Rolle's Theorem or the Mean Value Theorem.
Aull \cite{1967_Aull} proved quasi versions of both theorems for continuous functions, in which the secant slope is bounded from below and above by values of the symmetric derivative, and observed that they yield a Lipschitz-type estimate whenever the symmetric derivative is bounded; see \cite{1998_Sahoo_BOOK,2011_Sahoo} for further developments.

\begin{figure}[H]
  \centering
  \includegraphics[scale=1]{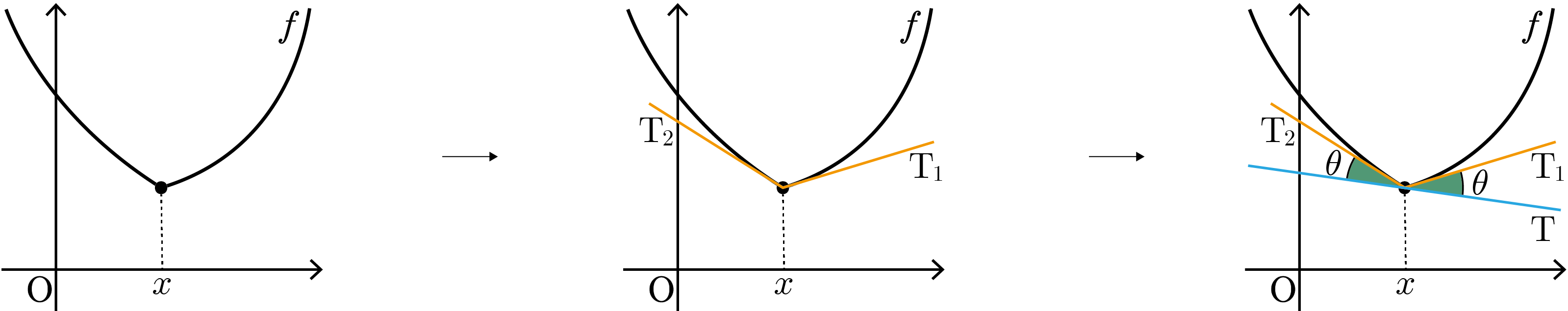}
  \caption{Motivation for specular derivatives}
  \label{Fig : Motivation for specular derivatives}
\end{figure}

In this paper, we introduce and study the \emph{regular specular derivative}, a generalized derivative that assigns a single slope to every point at which both one-sided derivatives exist and are finite.
The idea is geometric.
Suppose that $f$ has distinct finite one-sided derivatives at $x$, as in Figure \ref{Fig : Motivation for specular derivatives}.
The two one-sided tangent lines $\mathrm{T}_1$ and $\mathrm{T}_2$ at $(x,f(x))$ have slopes $f'_+(x)$ and $f'_-(x)$, respectively.
Think of $\mathrm{T}_2$ as a light ray arriving from the left and of $\mathrm{T}_1$ as the ray leaving to the right.
The mirror that turns one ray into the other is the line $\mathrm{T}$ through $(x,f(x))$ making equal angles with $\mathrm{T}_1$ and $\mathrm{T}_2$.
We define the regular specular derivative of $f$ at $x$ as the slope of $\mathrm{T}$; the word \enquote{specular} refers to this mirror.
At a point of differentiability the two tangent rays coincide with the tangent line, and so does the mirror.
Thus the regular specular derivative extends the classical derivative.

Whereas the symmetric derivative is the arithmetic mean of the two one-sided derivatives whenever they exist, the regular specular derivative is the slope obtained by averaging their angles of inclination.
This angular average commutes with rotations of the graph and obeys the reciprocal rule for inverse functions, two properties that the arithmetic mean lacks.
At each corner, the regular specular derivative can be calculated directly from the two finite one-sided slopes.

Initial ideas and an early treatment of specular differentiation were presented in \cite{2023_Jung}.
An independent difference-quotient definition and the associated general theory were developed in \cite{2026a_Jung}.
The present paper provides a self-contained geometric and analytical treatment of the regular case, in which both one-sided derivatives exist as finite real numbers; in this case the geometric definition agrees with the limit definition of \cite{2026a_Jung}, see Theorem \ref{Thm:sd_criterion}.
We revise and extend the earlier work in \cite{2023_Jung}, refining its definitions and terminology in light of \cite{2026a_Jung}, and we establish new results that appear in neither work.

The main results are as follows.
In $\mathbb{R}$, we derive three formulas for the regular specular derivative and establish inverse function and rotation formulas.
We prove Quasi-Rolle's Theorem and the Quasi-Mean Value Theorem.
Applying the Quasi-Mean Value Theorem, we establish a derivative-limit theorem.
In particular, twice regular specular differentiability implies $C^1$ regularity.
We also establish the Fundamental Theorem of Calculus for regular specular derivatives.

In $\mathbb{R}^n$, we extend the construction through coordinate phototangents and define specular tangent hyperplanes, regular specular gradients, and regular specular directional derivatives.
We show that continuity of the regular specular gradient implies classical differentiability.
Regular specular differentiability alone does not ensure that a regular specular directional derivative equals the inner product of the regular specular gradient with the direction vector.
Finally, we characterize the uniqueness of specular tangent hyperplanes using coordinate phototangents and the relation between the regular specular gradient and specular normal vectors.

The rest of the paper is organized as follows.
Section 2 develops the theory of regular specular derivatives for single-variable functions, including their geometric interpretation and applications to regularity and integration.
Section 3 treats functions of several variables, including specular tangent hyperplanes, gradients, directional derivatives, and normal vectors.

\section{Specular derivatives for single-variable functions}

Throughout this section, let $I$ be an open interval in $\mathbb{R}$.

\subsection{Definitions and properties}

\begin{definition}
  \label{def:sd}
  Let $f : I \to \mathbb{R}$ be a function, and let $x \in I$.
  We say that $f$ is \emph{regularly specularly differentiable} at $x$ if the right- and left-hand derivatives $f'_+(x)$ and $f'_-(x)$ both exist as finite real numbers.
  In this case, we define the following:
  \begin{enumerate}[label=\upshape(\alph*)]
      \item We define the \emph{phototangent} of $f$ at $x$ to be the function $\operatorname{pht}f:\mathbb{R}\to \mathbb{R}$ by
      \begin{equation*}
        \operatorname{pht}f(y)=
        \begin{cases}
          f'_+(x)(y - x) + f(x) & \text{if } y \geq x,\\
          f'_-(x)(y - x) + f(x) & \text{if } y < x.
        \end{cases}
      \end{equation*}
      \item The (\emph{first-order}) \emph{regular specular derivative} of $f$ at $x$, denoted by $f^{\sd}(x)$, is defined as the slope of the line through the two distinct intersection points of the graph of $\operatorname{pht}f$ and the unit circle $\partial B(\mathbf{z},1)$, where $\mathbf{z} = (x, f(x))$.
      \item We define the \emph{specular tangent line} to the graph of $f$ at the point $\left( x, f(x) \right)$, denoted by $\operatorname{stg}f$, to be the line passing through the point $\left( x, f(x) \right)$ with slope $f^{\sd}(x)$.
  \end{enumerate}
\end{definition}

We say that $f$ is regularly specularly differentiable on $I$ if it is regularly specularly differentiable at every point of $I$.
On a closed interval $[a, b]$, we require regular specular differentiability on $(a,b)$ and the existence of the finite one-sided derivatives $f'_+(a)$ and $f'_-(b)$.
At the endpoints, we set $f^{\sd}(a):=f'_+(a)$ and $f^{\sd}(b):=f'_-(b)$.

Note that regular specular differentiability at a point implies continuity at that point.
In this paper, we focus primarily on regularly specularly differentiable functions.
The notation $f^{\sd}$ in this section always refers to the regular specular derivative.

\begin{proposition}
  \label{Prop:sd_iff_pht_cont}
  Suppose that $f:I\to\mathbb{R}$ is regularly specularly differentiable at $x$.
  Its phototangent is continuous, and its graph meets every circle centered at $(x,f(x))$ with radius $r>0$ at exactly two points.
  The line through these points has a finite slope independent of $r$.
\end{proposition}

\begin{proof}
  Write $\alpha = f'_+(x)$, $\beta = f'_-(x)$, and $\mathbf{z} = (x, f(x))$.
  Let $r > 0$.
  The circle $\partial B(\mathbf{z},r)$ is given by $(y-x)^2 + (t - f(x))^2 = r^2$.
  Solving this equation with the defining expressions for the phototangent on $[x,\infty)$ and $(-\infty,x]$, respectively, gives the horizontal coordinates of the two intersection points:
  \begin{equation*}
    a := x + \frac{r}{\sqrt{1+\alpha^2}} 
    \qquad\text{and}\qquad
    b := x - \frac{r}{\sqrt{1+\beta^2}}.
  \end{equation*}
  Notice that $b < x < a$.
  To prove that $\operatorname{pht}f$ is continuous at $x$, take $\delta := \min\{ a-x, x-b \} > 0$.
  If $y \in (x-\delta, x]$, then
  \begin{equation*}
    |\operatorname{pht}f(y)-\operatorname{pht}f(x)|
    \leq
    |\operatorname{pht}f(b)-\operatorname{pht}f(x)|
    =
    \frac{r|\beta|}{\sqrt{1+\beta^2}}
    <
    r.
  \end{equation*}
  Similarly, if $y\in(x,x+\delta)$, then
  \begin{equation*}
    |\operatorname{pht}f(y)-\operatorname{pht}f(x)|
    \leq
    |\operatorname{pht}f(a)-\operatorname{pht}f(x)|
    =
    \frac{r|\alpha|}{\sqrt{1+\alpha^2}}
    <
    r.
  \end{equation*}
  Since $r > 0$ is arbitrary, $\operatorname{pht}f$ is continuous at $x$.
  It is also continuous away from $x$, since its restriction to each of $(-\infty,x)$ and $(x,\infty)$ is affine.
  Therefore, $\operatorname{pht}f$ is continuous on $\mathbb{R}$.

  The two intersection points are
  \begin{equation} \label{eq:intersection_btw_ball_pht}
    P_r = \left(x-\frac{r}{\sqrt{1+\beta^2}}, f(x)-\frac{r\beta}{\sqrt{1+\beta^2}}\right) 
    \qquad\text{and}\qquad
    Q_r = \left(x+\frac{r}{\sqrt{1+\alpha^2}}, f(x)+\frac{r\alpha}{\sqrt{1+\alpha^2}}\right).
  \end{equation}
  Their horizontal coordinates are distinct, so the line through them has a finite slope.
  The common factor $r$ cancels from the quotient defining this slope, proving that it is independent of $r$.
\end{proof}

\begin{corollary}\label{Crl : Linearity of phototangents}
  If $f,g:I\to\mathbb{R}$ are regularly specularly differentiable at $x$, then so is $f+g$, and their phototangents at $x$ satisfy
  \begin{equation*}
    \operatorname{pht}(f+g)=\operatorname{pht}f+\operatorname{pht}g.
  \end{equation*}
\end{corollary}

\begin{proof}
  This follows from $(f+g)'_\pm(x)=f'_\pm(x)+g'_\pm(x)$ and the definition of the phototangent.
\end{proof}

\begin{example}
  The function $f(x)=|x|$ is regularly specularly differentiable on $\mathbb{R}$, with $f^{\sd}(x)=\operatorname{sgn}(x)$.
  Thus a regular specular derivative need not be continuous.
  The sign function itself is not regularly specularly differentiable at $0$, and its phototangent at $0$ is not defined in the present framework.
  Likewise, the function $g$ given by $g(x)=|x|$ for $x \neq 0$ and $g(0)=1$ is not regularly specularly differentiable at $0$, since it is not continuous there.
\end{example}

\begin{remark}
  \label{rmk:properties_stl}
  For a function regularly specularly differentiable at $x$, its specular tangent line is given by
  \begin{equation*}
    \operatorname{stg}f(y)=f(x)+f^{\sd}(x)(y-x),\qquad y\in\mathbb{R}.
  \end{equation*}
  In particular, $\operatorname{stg}f(x)=f(x)$ and $(\operatorname{stg}f)'(x)=f^{\sd}(x)$.
\end{remark}

\begin{figure}[H]
  \centering
  \includegraphics[scale=1]{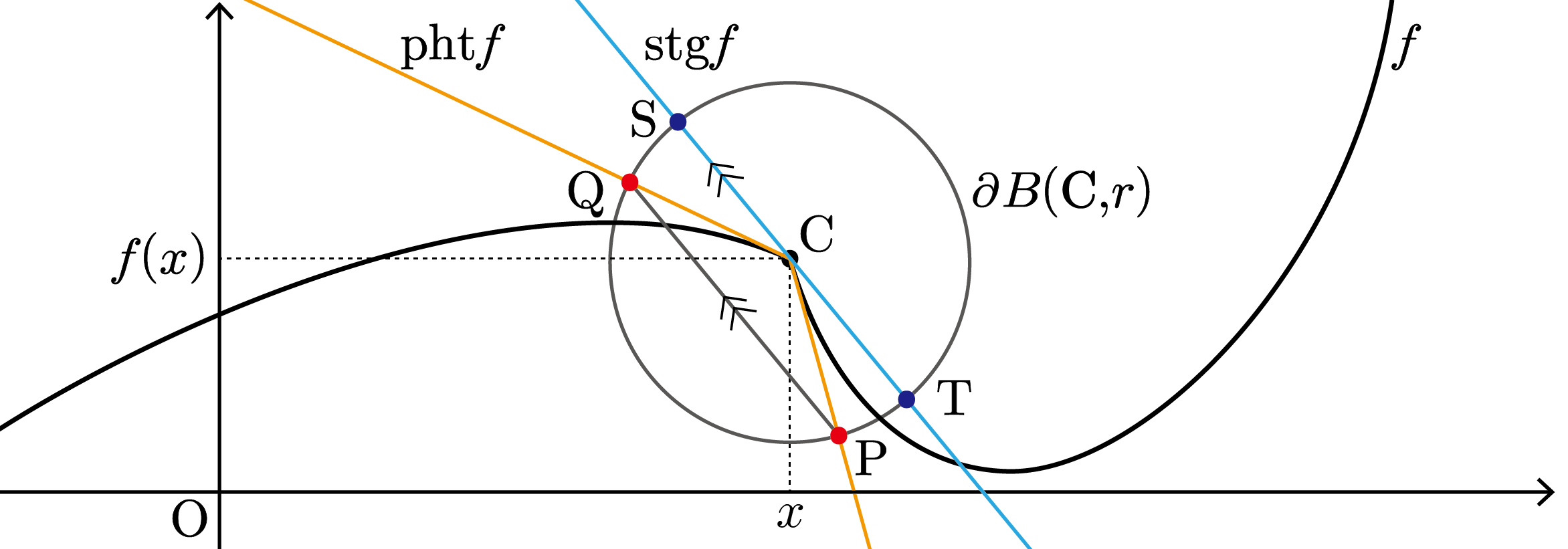}
  \caption{Basic concepts concerning specular derivatives}
  \label{Fig : basic concepts concerning specular derivatives}
\end{figure}

In Figure \ref{Fig : basic concepts concerning specular derivatives}, the function $f$ is regularly specularly differentiable at $x$ but is not classically differentiable there.
Let $\operatorname{pht}f$ be the phototangent of $f$ at $x$, and let $\operatorname{stg}f$ be the specular tangent line to the graph of $f$ at $(x ,f(x))$.
By Proposition \ref{Prop:sd_iff_pht_cont}, the phototangent is continuous at $x$.
Imagine that you shoot a light ray toward a mirror.
The words \enquote{specular} in specular tangent line and \enquote{photo} in phototangent stand for the mirror $\operatorname{stg}f$ and the light ray $\operatorname{pht}f$, respectively.
Write $C = (x, f(x))$, and adopt the notation in Figure \ref{Fig : basic concepts concerning specular derivatives}.
Since $CP = CQ$, the triangle $CPQ$ is isosceles.
The line $ST$ passes through $C$ and is parallel to $PQ$, and hence
\begin{equation*}
  \angle CPQ=\angle CQP=\angle SCQ=\angle TCP.
\end{equation*}
Thus the specular tangent line represents a mirror for which the incident and reflected light rays make equal angles with the mirror.

Moreover, Figure \ref{Fig : basic concepts concerning specular derivatives} illustrates why the specular derivative is independent of $r$, as proved in Proposition \ref{Prop:sd_iff_pht_cont}.
For $0 < r < 1$, let $U = (u_1,u_2)$ and $V = (v_1,v_2)$ denote the intersection points of $\operatorname{pht}f$ and $\partial B(C,1)$, where $u_1 < x < v_1$.
Since the ratios $\overline{CU}:\overline{CQ}=1:r$ and $\overline{CV}:\overline{CP} = 1:r$ hold, the line segments $\overline{PQ}$ and $\overline{UV}$ are parallel.
The case when $r > 1$ can be treated similarly.

We suggest three ways to calculate a regular specular derivative.
The first formula expresses it as a limit of slopes determined by secants to the graph.
Following the notation in \cite{2026a_Jung}, let the auxiliary function $\mathcal{A}: \mathbb{R} \times \mathbb{R} \times (0, \infty) \to \mathbb{R}$ be given by
\begin{equation}
  \label{eq:auxiliary_A}
  \mathcal{A}(a,b,c):= \frac{a\sqrt{b^2+c^2}+b\sqrt{a^2+c^2}}{c\sqrt{a^2+c^2}+c\sqrt{b^2+c^2}}.
\end{equation}

We derive the formula below directly from the geometry of the phototangent.

\begin{theorem}
  \label{Thm:sd_criterion}
  Let $f : I \to \mathbb{R}$ be a function, and let $x\in I$.
  If $f$ is regularly specularly differentiable at $x$, then
  \begin{equation*}
    f^{\sd}(x) = \lim_{h\searrow0} \mathcal{A}\bigl(f(x+h) - f(x), f(x) - f(x-h),h\bigr).
  \end{equation*}
\end{theorem}

\begin{figure}[H]
  \centering
  \includegraphics[scale=1]{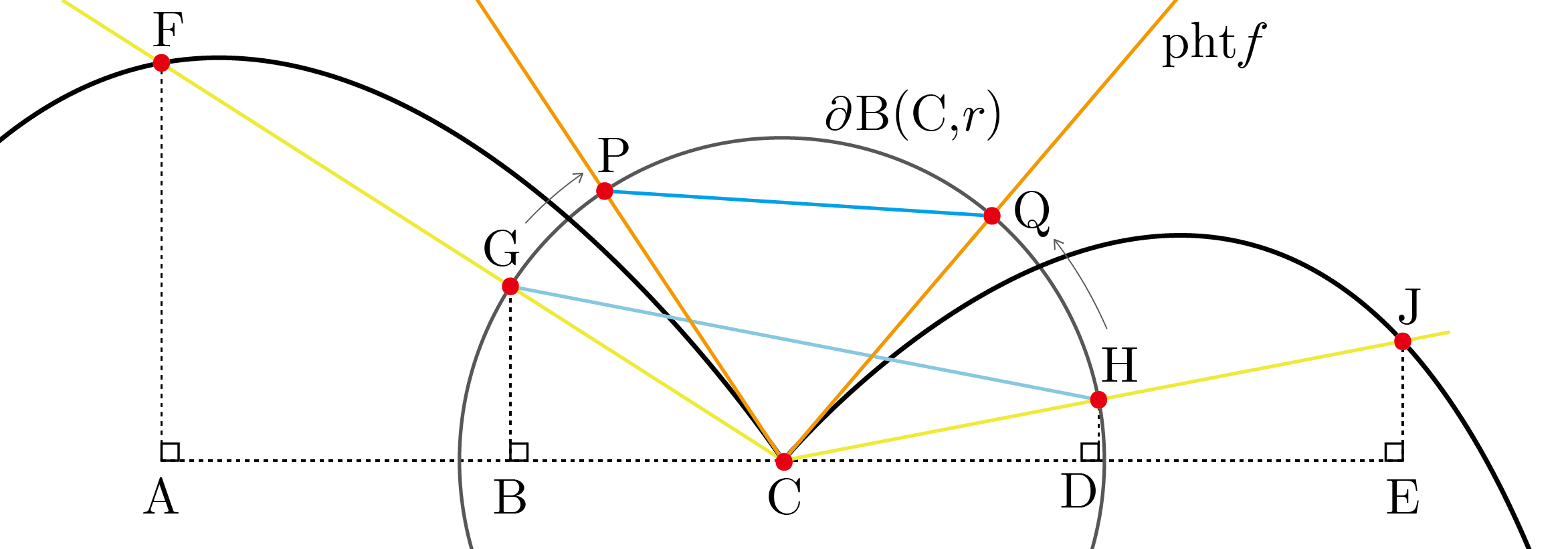}
  \caption{The slope of the line GH converges to the specular derivative of $f$ at $x$}
  \label{Fig : The slope of the line GH converges the specular derivative of f at x}
\end{figure}

\begin{proof}
  Let $h>0$ be sufficiently small that $x-h, x+h \in I$.
  Denote $A:=(x-h,f(x))$, $C:=(x,f(x))$, $E:=(x+h,f(x))$, $F:=(x-h,f(x-h))$, and $J:=(x+h,f(x+h))$.
  Choose $r=r(h)$ such that $0<r(h)<h$.
  By Proposition  \ref{Prop:sd_iff_pht_cont}, the graph of $\operatorname{pht}f$ intersects the circle $\partial B(C, r)$  at exactly two points, denoted by $P := (p_1, p_2)$ and $Q := (q_1, q_2)$, where $p_1 < x < q_1$ .
  Note that the slopes of $\overline{CP}$ and $\overline{CQ}$ are equal to $f'_-(x)$ and $f'_+(x)$, respectively.
  Let $G:=(g_1,g_2)$ and $H:=(h_1,h_2)$ be the intersection points of $\partial B(C,r)$ with the line segments $\overline{CF}$ and $\overline{CJ}$, respectively.
  Since $0<r<h$, these points exist and satisfy $g_1 < x < h_1$.
  Denote $B:=(g_1,f(x))$ and $D:=(h_1,f(x))$.
  See Figure \ref{Fig : The slope of the line GH converges the specular derivative of f at x}.

  Observe the following pairs of similar right triangles:
  \begin{equation*}
    \triangle CAF\sim\triangle CBG \qquad\text{and}\qquad \triangle CEJ\sim\triangle CDH.
  \end{equation*}
  By the basic properties of similar right triangles, we obtain
  \begin{equation*}
    g_1 = x - \frac{rh}{\sqrt{(f(x-h)-f(x))^2+h^2}}
    \qquad\text{and}\qquad
    h_1 = x+\frac{rh}{\sqrt{(f(x+h)-f(x))^2+h^2}}.
  \end{equation*}
  Similarly, one can find that
  \begin{equation*}
    g_2 = f(x) + \frac{r(f(x-h)-f(x))}{\sqrt{(f(x-h)-f(x))^2+h^2}} 
    \qquad\text{and}\qquad
    h_2 = f(x)+\frac{r(f(x+h)-f(x))}{\sqrt{(f(x+h)-f(x))^2+h^2}}.
  \end{equation*}
  These formulas also hold directly when either secant is horizontal, in which case the corresponding triangles are degenerate.
  The slope of $\overline{GH}$ is given by
  \begin{equation*}
    \frac{h_2-g_2}{h_1-g_1} 
    =
    \frac{ \displaystyle\frac{f(x+h)-f(x)}{\sqrt{(f(x+h)-f(x))^2+h^2}} -\displaystyle\frac{f(x-h)-f(x)}{\sqrt{(f(x-h)-f(x))^2+h^2}} }{ \displaystyle\frac{h}{\sqrt{(f(x+h)-f(x))^2+h^2}} +\displaystyle\frac{h}{\sqrt{(f(x-h)-f(x))^2+h^2}} },
  \end{equation*}
  which is independent of $r$.
  Multiplying the numerator and denominator by
  \begin{equation*}
    \sqrt{(f(x+h)-f(x))^2+h^2}\cdot \sqrt{(f(x-h)-f(x))^2+h^2}
  \end{equation*}
  implies
  \begin{align*}
    \frac{h_2-g_2}{h_1-g_1}
    &=\frac{(f(x+h)-f(x))\sqrt{(f(x-h)-f(x))^2+h^2} - (f(x-h)-f(x))\sqrt{(f(x+h)-f(x))^2+h^2}}{h\sqrt{(f(x+h)-f(x))^2+h^2} + h\sqrt{(f(x-h)-f(x))^2+h^2}} \\
    &=\mathcal{A}\bigl(f(x+h)-f(x),f(x)-f(x-h),h\bigr).
  \end{align*}
  The slopes of $\overline{CF}$ and $\overline{CJ}$ converge to the finite slopes $f'_-(x)$ and $f'_+(x)$ of $\overline{CP}$ and $\overline{CQ}$, respectively.
  Consequently, the normalized vectors $\frac{1}{r}(G-C)$ and $\frac{1}{r}(H-C)$ converge to $\frac{1}{r}(P-C)$ and $\frac{1}{r}(Q-C)$, whose expressions are independent of $r$.
  In particular,
  \begin{equation*}
    \frac{h_1-g_1}{r} \to \frac{1}{\sqrt{1+f'_-(x)^2}}+ \frac{1}{\sqrt{1+f'_+(x)^2}}>0.
  \end{equation*}
  It follows that the slope of $\overline{GH}$ converges to that of $\overline{PQ}$, which equals the regular specular derivative $f^{\sd}(x)$, as $h \searrow 0$ .
\end{proof}

Thus, in the regular case, the geometric definition in Definition \ref{def:sd} agrees with the limit definition of the specular derivative in \cite{2026a_Jung}.

The symmetry of the expression in Theorem \ref{Thm:sd_criterion} gives the following consequence.

\begin{corollary}
  Let $f : I\to \mathbb{R}$ be a function, and $x_0 \in I$ be a point.
  Assume $f$ is regularly specularly differentiable at $x_0$.
  If $f$ is symmetric about $x=x_0$ in a neighborhood of $x_0$, that is, there exists $\delta > 0$ such that 
  \begin{equation*} 
    f(x_0 - x) = f(x_0 + x)
  \end{equation*}
  for all $x\in(0,\delta)$, then $f^{\sd} ( x_0 ) = 0$.
\end{corollary}

\begin{proof}
  The numerator in \eqref{eq:auxiliary_A} vanishes whenever $a = -b$.
  Apply Theorem \ref{Thm:sd_criterion} with $f(x_0+x) = f(x_0-x)$.
\end{proof}

\begin{example}
  For the ReLU function $f(x)=\frac{1}{2}\left(x+|x|\right)$, one can calculate $f^{\sd}(0)=-1+\sqrt{2}$ by using Theorem \ref{Thm:sd_criterion}.
\end{example}

In order to calculate regular specular derivatives more conveniently, we suggest the second formula using ordinary one-sided derivatives.

\begin{proposition} \label{Prop : Calculating spd}
  Suppose that $f : I \to \mathbb{R}$ is regularly specularly differentiable at a point $x \in I$.
  Write $f'_+(x)=:\alpha$ and $f'_-(x)=:\beta$.
  Then we have
  \begin{equation} \label{Prop : Calculating spd formula}
    f^{\sd} (x) =
    \begin{cases} 
      \displaystyle \frac{\alpha\beta-1 + \sqrt{\left(\alpha^2 + 1\right)\left( \beta^2 + 1 \right)}}{\alpha+\beta} & \text{if } \alpha+\beta\neq 0,\\[1em] 
      0 & \text{if } \alpha+\beta=0. 
    \end{cases}      
  \end{equation}
\end{proposition}

\begin{proof}
  By \eqref{eq:intersection_btw_ball_pht}, the slope of the line through the two circle intersections is 
  \begin{equation}
    \label{eq:regular-weighted-formula}
    f^{\sd}(x) = \frac{\alpha\sqrt{1+\beta^2}+\beta\sqrt{1+\alpha^2}}{\sqrt{1+\alpha^2}+\sqrt{1+\beta^2}}.
  \end{equation}
  If $\alpha+\beta=0$, the numerator is zero.
  Otherwise, write $u=\sqrt{1+\alpha^2}$ and $v=\sqrt{1+\beta^2}$.
  Direct expansion gives
  \begin{equation*}
    (\alpha+\beta)(\alpha v+\beta u)=(\alpha\beta-1+uv)(u+v).
  \end{equation*}
  Dividing by $(\alpha+\beta)(u+v)$ yields \eqref{Prop : Calculating spd formula}.
\end{proof}

The following corollary records some consequences of this formula.

\begin{corollary} \label{cor:calculating_sd}
  If $f: I \to \mathbb{R}$ is regularly specularly differentiable at $x \in I$, then the following statements hold:
  \begin{enumerate}[label=\upshape(\alph*)]
    \item \label{cor:calculating_sd-1} 
    
    $f^{\sd}(x) = 0$ if and only if $f'_+(x) + f'_-(x) = 0$. 
    \item \label{cor:calculating_sd-2} 
    
    Signs of $f^{\sd}(x)$ and $f'_+(x) + f'_-(x)$ are equal, i.e., $\operatorname{sgn}\left( f^{\sd}(x) \right) = \operatorname{sgn}\left( f'_+(x) + f'_-(x) \right)$.
    \item \label{cor:calculating_sd-3} 
    
    $\displaystyle |f^{\sd}(x)| \leq \frac{|f'_+(x) + f'_-(x)|}{2}$.

    \item \label{cor:calculating_sd-4}
    If $f$ is Lipschitz with constant $L$ on a neighborhood of $x$, then $|f^{\sd}(x)|\leq L$.
  \end{enumerate}
\end{corollary}

\begin{proof}
  Write $\alpha = f'_+(x)$ and $\beta = f'_-(x)$.
  If $\alpha+\beta=0$, Proposition \ref{Prop : Calculating spd} gives $f^{\sd}(x)=0$.
  If $\alpha+\beta\ne0$, then
  \begin{equation*}
    (1+\alpha^2)(1+\beta^2) - (1-\alpha\beta)^2
    =
    (\alpha+\beta)^2
    >
    0,
  \end{equation*}
  and hence $\alpha\beta-1+\sqrt{(1+\alpha^2)(1+\beta^2)} > 0$.
  The sign of \eqref{Prop : Calculating spd formula} is therefore the sign of $\alpha+\beta$.
  This proves \ref{cor:calculating_sd-1} and \ref{cor:calculating_sd-2}.

  For \ref{cor:calculating_sd-3}, if $\alpha+\beta=0$, then both sides of the inequality are zero.
  Now suppose that $\alpha+\beta\neq0$.
  By the positivity of the numerator established above and \eqref{Prop : Calculating spd formula}, we have
  \begin{equation*}
    |f^{\sd}(x)|
    =
    \frac{\alpha\beta-1+\sqrt{(\alpha^2+1)(\beta^2+1)}}{|\alpha+\beta|}.
  \end{equation*}
  The Arithmetic Mean--Geometric Mean Inequality implies that
  \begin{equation*}
    \alpha\beta-1+\sqrt{(\alpha^2+1)(\beta^2+1)}
    \leq
    \alpha\beta-1+\frac{(\alpha^2+1)+(\beta^2+1)}{2}
    =
    \frac{(\alpha+\beta)^2}{2}.
  \end{equation*}
  Dividing by $|\alpha+\beta|$ proves \ref{cor:calculating_sd-3}.

  To prove \ref{cor:calculating_sd-4}, suppose that $f$ is Lipschitz with constant $L$ on a neighborhood of $x$.
  Choose $r>0$ such that $(x-r,x+r)\subset I$ and the Lipschitz condition holds on this interval.
  For $0<|h|<r$, we have
  \begin{equation*}
    \left|\frac{f(x+h)-f(x)}{h}\right|
    \leq
    \frac{L|h|}{|h|}
    =
    L.
  \end{equation*}
  Taking the right- and left-hand limits yields $|\alpha|\leq L$ and $|\beta|\leq L$.
  By \ref{cor:calculating_sd-3} and the triangle inequality,
  \begin{equation*}
    |f^{\sd}(x)|
    \leq
    \frac{|\alpha+\beta|}{2}
    \leq
    \frac{|\alpha|+|\beta|}{2}
    \leq
    L.
  \end{equation*}
\end{proof}

\begin{remark} \label{Rmk : Specular derivatives may not be linear}
  Regular specular derivatives are not linear. 
  For instance, consider the ReLU function $f(x) =\frac{1}{2}\left( x + |x| \right)$ for $x \in \mathbb{R}$.
  First, we find that 
  \begin{equation*} 
    \left( 2 f \right)^{\sd}(0) = \frac{-1+\sqrt{5}}{2} \neq 2 \left( -1 + \sqrt{2} \right) = 2 f^{\sd}(0)
  \end{equation*}
  by using \eqref{Prop : Calculating spd formula} with $(2f)'_+(0) = 2$, $(2f)'_-(0) = 0$, $f'_+(0) = 1$, and $f'_-(0) = 0$.

  Also, take the smooth function $g(x)=x$ for $x\in\mathbb{R}$.
  Second, one can calculate that
  \begin{equation*}
    f^{\sd}(0)+g^{\sd}(0)
    =
    \bigl(-1+\sqrt{2}\bigr) + 1
    =
    \sqrt{2} 
    \neq
    \frac{1+\sqrt{10}}{3}
    =
    (f+g)^{\sd}(0).
  \end{equation*}
  Furthermore, regular specular derivatives need not obey the classical chain rule.
  Consider the composite function $f\circ h$ with $h(x)=2x$.
  Writing $y=h(x)$, we have $x = 0$  if and only if $y = 0$.
  Then we find that
  \begin{equation*}
    (f\circ h)^{\sd}(0)
    =
    \frac{-1+\sqrt{5}}{2}
    \neq
    (-1+\sqrt{2})2
    =
    f^{\sd}(h(0))h^{\sd}(0).
  \end{equation*}
\end{remark}

As stated in the next theorem, regular specular derivatives generalize classical derivatives.

\begin{theorem} \label{Thm : ordinary dervatives and specular derivatives}
  Let $f: I \to \mathbb{R}$ be a function, and let $x \in I$ be a point.
  \begin{enumerate}[label=\upshape(\alph*)]
    \item
    \label{Item : differentiability implies specularly differentiability}
    If $f$ is differentiable at $x$, then $f$ is regularly specularly differentiable at $x$ and $f'(x)=f^{\sd} (x)$. 
    
    \item
    \label{Item : differentiability iff pht differentiability}
    If $f$ is regularly specularly differentiable at $x$, then $f$ is differentiable at $x$ if and only if its phototangent is differentiable at $x$. 
  \end{enumerate}
\end{theorem}

\begin{proof}
  To prove \ref{Item : differentiability implies specularly differentiability}, assume $f$ is differentiable at $x$.
  Then $f'_+(x)=f'_-(x)=f'(x)\in\mathbb{R}$, and hence $f$ is regularly specularly differentiable at $x$ by definition.
  On the one hand, if $f'_+(x)=f'_-(x)=0$, we see that $f'(x) = f'_+(x) = 0 = f^{\sd}(x)$, using Proposition \ref{Prop : Calculating spd}.
  On the other hand, if $f'_+(x)\neq 0$ and $f'_-(x) \neq 0$, then $f'_+ (x) + f'_- (x) = 2 f'_+ (x) \neq 0$.
  Writing $\alpha := f'_+ (x)$, one can calculate
  \begin{equation*}
    f'(x) = f'_+(x) = \frac{\alpha^2 - 1 + \sqrt{\left( \alpha^2 + 1 \right)^2}}{2 \alpha} = f^{\sd}(x)
  \end{equation*}
  by applying Proposition \ref{Prop : Calculating spd} again.
  This proves \ref{Item : differentiability implies specularly differentiability}.

  Now, assume $f$ is regularly specularly differentiable at $x$, and let $\operatorname{pht}f$ be its phototangent to show \ref{Item : differentiability iff pht differentiability}.
  First, suppose $f$ is differentiable at $x$.
  Since $f'_+(x)=f'_-(x)$, we conclude that $\operatorname{pht}f$ is a polynomial of degree $1$ or less.
  Conversely, assume the phototangent $\operatorname{pht}f$ is differentiable at $x$.
  Observing that $\left(\operatorname{pht}f\right)'_+(x)=\left(\operatorname{pht}f\right)'_-(x)$, we have $f'_+(x)=f'_-(x)$, which implies that $f$ is differentiable at $x$.
\end{proof}

We suggest the third formula calculating specular derivatives and its applications.
Following the notation in \cite{2026a_Jung}, let $\mathcal{B}:\mathbb{R}\times\mathbb{R}\to\mathbb{R}$ be given by
\begin{equation*}
  \mathcal{B}(\alpha, \beta) := \tan \left( \frac{\arctan \alpha + \arctan \beta}{2} \right).
\end{equation*}

\begin{lemma} \label{lem:sd_tan_arctan} 
  If $f:I \to \mathbb{R}$ is regularly specularly differentiable at a point $x \in I$, then 
  \begin{equation*}
    f^{\sd}(x)
    =
    \mathcal{B} \bigl(f'_+(x), f'_-(x)\bigr).
  \end{equation*}
\end{lemma}

\begin{proof}
  Write $\theta_1=\arctan f'_+(x)$ and $\theta_2=\arctan f'_-(x)$.
  Then $\theta_1,\theta_2\in\left(-\frac{\pi}{2},\frac{\pi}{2}\right)$.
  By \eqref{eq:intersection_btw_ball_pht}, the chord joining the two intersections of the phototangent with the unit circle centered at $(x,f(x))$ has slope
  \begin{equation*}
    f^{\sd}(x)
    =
    \frac{\sin\theta_1+\sin\theta_2}{\cos\theta_1+\cos\theta_2}
    =
    \tan\left(\frac{\theta_1+\theta_2}{2}\right)
    =
    \mathcal{B}\bigl(f'_+(x),f'_-(x)\bigr).
  \end{equation*}
  The denominator is positive because both angles belong to $\left(-\frac{\pi}{2},\frac{\pi}{2}\right)$.
  The second equality follows from the sum-to-product identities, and the last equality follows from the definition of $\mathcal{B}$.
\end{proof}

Reflection of a graph across the line $y=x$ gives the graph of its inverse.
For regular specular derivatives, this geometric relation yields the reciprocal rule under the following assumptions.

\begin{theorem}[Inverse function rule]
  \label{thm:regular-inverse}
  Let $f:I\to\mathbb{R}$ be continuous and strictly monotone, and let $x\in I$.
  Suppose that $f$ is regularly specularly differentiable at $x$ and that $f'_+(x)f'_-(x) \neq 0$.
  Then $f^{-1}:f(I)\to I$ is regularly specularly differentiable at $f(x)$, and
  \begin{equation}
    \label{eq:regular-inverse}
    \bigl(f^{-1}\bigr)^{\sd}(f(x))=\frac{1}{f^{\sd}(x)}.
  \end{equation}
\end{theorem}

\begin{proof}
  Write $y=f(x)$, $g=f^{-1}$, $\alpha=f'_+(x)$, and $\beta=f'_-(x)$.
  Strict monotonicity and the nonzero assumptions imply that $\alpha$ and $\beta$ have the same sign.
  In particular, $f^{\sd}(x) \neq 0$ by Corollary \ref{cor:calculating_sd}.
  The image $f(I)$ is an open interval, and $g$ is continuous there.
  For $z\in f(I)\setminus\{y\}$, set $t=g(z)$.
  Then $t\in I\setminus\{x\}$ and $z=f(t)$, and hence
  \begin{equation*}
    \frac{g(z)-g(y)}{z-y}
    =
    \frac{t-x}{f(t)-f(x)}
    =
    \left(\frac{f(t)-f(x)}{t-x}\right)^{-1}.
  \end{equation*}
  If $f$ is increasing, $z>y$ if and only if $t>x$.
  If $f$ is decreasing, $z>y$ if and only if $t<x$.
  Taking one-sided limits therefore gives
  \begin{equation*}
    \bigl(g'_+(y),g'_-(y)\bigr)
    =
    \begin{cases}
      \displaystyle \left(\frac{1}{\alpha},\frac{1}{\beta}\right) & \text{if } f \text{ is increasing},\\[1em]
      \displaystyle \left(\frac{1}{\beta},\frac{1}{\alpha}\right) & \text{if } f \text{ is decreasing}.
    \end{cases}
  \end{equation*}
  These derivatives are finite, and hence $g$ is regularly specularly differentiable at $y$.

  For $t\neq0$, we have
  \begin{equation*}
    \arctan\frac{1}{t}
    =\operatorname{sgn}(t)\frac{\pi}{2}-\arctan t.
  \end{equation*}
  Since $\alpha$ and $\beta$ have the same sign, the definition of $\mathcal{B}$ yields
  \begin{equation*}
    \mathcal{B}\left(\frac{1}{\alpha},\frac{1}{\beta}\right)
    =
    \tan\left( \operatorname{sgn}(\alpha)\frac{\pi}{2} - \frac{\arctan\alpha+\arctan\beta}{2} \right)
    =
    \frac{1}{\mathcal{B}(\alpha,\beta)}.
  \end{equation*}
  By Lemma \ref{lem:sd_tan_arctan} and the symmetry of $\mathcal{B}$, we conclude that
  \begin{equation*}
    g^{\sd}(y)
    =
    \mathcal{B}\bigl(g'_+(y),g'_-(y)\bigr)
    =
    \mathcal{B}\left(\frac{1}{\alpha},\frac{1}{\beta}\right)
    =
    \frac{1}{\mathcal{B}(\alpha,\beta)}
    =
    \frac{1}{f^{\sd}(x)}.
  \end{equation*}
\end{proof}

\begin{remark}
  \label{rmk:inverse-nonzero-slopes}
  The condition $f^{\sd}(x)\neq 0$ in Theorem \ref{thm:regular-inverse} alone does not ensure regular specular differentiability of the inverse.
  For example, the strictly increasing function
  \begin{equation*}
    f(t)=
    \begin{cases}
      t & \text{if } t\le0,\\
      t^2 & \text{if } t>0
    \end{cases}
  \end{equation*}
  satisfies $f'_+(0)=0$, $f'_-(0)=1$, and $f^{\sd}(0)=\sqrt{2}-1\neq0$.
  However, $f^{-1}(s)=\sqrt{s}$ for $s>0$, so its right-hand difference quotient at $0$ tends to $+\infty$.
\end{remark}

The following proposition shows that the specular tangent direction rotates through the same angle as the two one-sided tangent directions.

\begin{proposition}[Rotation formula]
  Let $f:I\to\mathbb{R}$ be regularly specularly differentiable at a point $x\in I$.
  Fix $\theta\in\mathbb{R}$ and define
  \begin{equation*}
    R(t):=\frac{\sin\theta+t\cos\theta}{\cos\theta-t\sin\theta}
  \end{equation*}
  whenever the denominator is nonzero.
  If $\cos\theta - f'_+(x) \sin\theta > 0$ and $\cos\theta - f'_-(x) \sin\theta > 0$, then 
  \begin{equation*}
    \mathcal{B}\left(R\bigl(f'_+(x)\bigr), R\bigl(f'_-(x)\bigr) \right)
    =
    R\left(\mathcal{B} \bigl( f'_+(x), f'_-(x) \bigr) \right)
    =
    R\bigl(f^{\sd}(x)\bigr).
  \end{equation*}
\end{proposition}

\begin{proof}
  Write $\alpha = f'_+(x)$, $\beta = f'_-(x)$, $\theta_1=\arctan\alpha$, and $\theta_2=\arctan\beta$.
  By periodicity, we may assume that $-\pi\le\theta<\pi$.
  The assumptions yield
  \begin{equation*}
    \cos(\theta_1+\theta)
    =
    \frac{\cos\theta-\alpha\sin\theta}{\sqrt{1+\alpha^2}} > 0 
    \qquad\text{and}\qquad
    \cos(\theta_2+\theta)
    =
    \frac{\cos\theta-\beta\sin\theta}{\sqrt{1+\beta^2}} > 0.
  \end{equation*}
  Since $\theta_1,\theta_2\in(-\frac{\pi}{2},\frac{\pi}{2})$, it follows that $-\frac{\pi}{2}<\theta_1+\theta<\frac{\pi}{2}$ and $-\frac{\pi}{2}<\theta_2+\theta<\frac{\pi}{2}$.
  Thus the tangent addition formula yields $\arctan R(\alpha) = \theta_1 + \theta$ and $\arctan R(\beta) = \theta_2 + \theta$.
  Consequently,
  \begin{equation*}
    \mathcal{B}\bigl(R(\alpha),R(\beta)\bigr)
    = \tan\left( \frac{(\theta_1+\theta)+(\theta_2+\theta)}{2} \right) 
    = \tan\left(\frac{\theta_1+\theta_2}{2}+\theta\right)
    = R\left(\tan \left( \frac{\theta_1+\theta_2}{2} \right)\right)
    = R\bigl(\mathcal{B}(\alpha,\beta)\bigr).
  \end{equation*}
  The last desired equality follows from Lemma~\ref{lem:sd_tan_arctan}.
\end{proof}

\subsection{The Quasi-Mean Value Theorem}

A regularly specularly differentiable function is classically differentiable Lebesgue almost everywhere by the Denjoy--Young--Saks theorem.
We seek conditions under which classical differentiability holds everywhere.
To this end, we first establish a weakened version of the Mean Value Theorem for regular specular derivatives.

\begin{example}
  Specular derivatives do not satisfy the classical Rolle's Theorem or the classical Mean Value Theorem.
  For the Mean Value Theorem, take $f:[-1,1]\to\mathbb{R}$ defined by $f(x)=x+|x|$.
  Its secant slope on $[-1,1]$ is $1$, whereas $f^{\sd}$ takes only the values $0$, $2$, and $\frac{1}{2}(\sqrt{5}-1)$.
  For Rolle's Theorem, consider $g:[0,3]\to\mathbb{R}$ defined by $g(x)=2x$ for $0\le x\le1$ and $g(x)=3-x$ for $1<x\le3$.
  Then $g(0)=g(3)=0$, but $g^{\sd}$ takes only the nonzero values $2$, $-1$, and $\sqrt{10}-3$.
\end{example}

\begin{lemma} \label{Lem : continuous and existence}
  Let $f$ be a continuous function on $[a, b] \subset \mathbb{R}$.
  Assume $f$ is regularly specularly differentiable in $(a,b)$.
  Then the following properties hold:
  \begin{enumerate}[label=\upshape(\alph*)]
    \item If $f(a)<f(b)$, then there exists $c_1 \in (a,b)$ such that $f^{\sd}(c_1) \geq 0$.
    
    \item 
    \label{Lem : continuous and existence (b)}
    If $f(a)>f(b)$, then there exists $c_2 \in (a,b)$ such that $f^{\sd}(c_2) \leq 0$. 
  \end{enumerate}
\end{lemma}

\begin{proof}
  First of all, assume $f(a)<f(b)$.
  Throughout the proof, $k$ denotes a real number with $f(a)<k<f(b)$.
  Since $f$ is continuous and $f(a) < k < f(b)$, by the intermediate value theorem, there exists $x \in (a, b)$ such that $f(x) = k$, and hence the set
  \begin{equation*}
    K := \{ x \in [a, b] : f(x) > k \}
  \end{equation*}
  contains $b$ and is therefore nonempty.
  Since $K$ is nonempty and bounded below, define $c_1:=\inf K$.
  Continuity at the endpoints gives $a<c_1<b$, and continuity at $c_1$ gives $f(c_1)=k$.
  Then, by the definition of $c_1 = \inf K$, there exists a sequence $\{ h_j \}_{j = 1}^\infty$ such that $h_j > 0$ and $c_1 + h_j \in K$ for all $j \in \mathbb{N}$, and $h_j \searrow 0$ as $j \to \infty$.
  For each $j \in \mathbb{N}$, it holds that
  \begin{equation*}
    \frac{f(c_1 + h_j) - f(c_1)}{h_j} > \frac{k - f(c_1)}{h_j} \geq 0
  \end{equation*}
  since $f(c_1) \leq k$.
  Since the finite right-hand derivative exists, taking the limit as $j \to \infty$ yields
  \begin{equation*}
    f'_+(c_1) = \lim_{j \to \infty} \frac{f(c_1 + h_j) - f(c_1)}{h_j} \geq 0.
  \end{equation*}
  Hence, $f'_+(c_1) \geq 0$.
  Since $f(c_1 - h) \leq k$ for any sufficiently small $h > 0$, we have
  \begin{equation*}
    \frac{f(c_1) - f(c_1 - h)}{h} \geq \frac{f(c_1) - k}{h} \geq 0.
  \end{equation*}
  Taking the limit as $h \searrow 0$ yields
  \begin{equation*}
    f'_-(c_1) = \lim_{h \searrow 0} \frac{f(c_1)- f(c_1 - h)}{h} \geq 0.
  \end{equation*}

  On the one hand, assume $f'_+(c_1) + f'_-(c_1)=0$.
  Then $f^{\sd}(c_1)=0$ due to Proposition \ref{Prop : Calculating spd}.
  On the other hand, suppose $f'_+(c_1) + f'_-(c_1) \neq 0$.
  One can estimate that
  \begin{equation*}
    f^{\sd}(c_1) \geq \frac{f'_+(c_1)f'_-(c_1)}{f'_+(c_1)+f'_-(c_1)} \geq 0,
  \end{equation*}
  using Proposition \ref{Prop : Calculating spd}.
  Hence, we conclude that $f^{\sd}(c_1) \geq 0$.

  Similarly, the proof of the reversed inequalities in \ref{Lem : continuous and existence (b)} can be shown.
\end{proof}

\begin{theorem}[Quasi-Rolle's Theorem]
  \label{thm:Quasi-Rolle}
  Let $f : [a, b] \to \mathbb{R}$ be a continuous function on $[a,b]$. 
  Suppose $f$ is regularly specularly differentiable in $(a,b)$ and $f(a) = f(b) = 0$.
  Then there exist $c_1$, $c_2 \in (a,b)$ such that $f^{\sd}(c_2)\leq  0 \leq  f^{\sd}(c_1)$.
\end{theorem}

\begin{proof} 
  If $f \equiv 0$, the conclusion follows.
  Now, suppose $f \not\equiv 0$.
  The hypothesis implies three cases; there exists either $a^{\ast} \in (a,b)$ such that $f\left(a^{\ast}\right)>0$, or $b^{\ast} \in (a,b)$ such that $f \left(b^{\ast} \right) <0$, or both.
  If such $a^{\ast}$ exists, using Lemma \ref{Lem : continuous and existence} on $[a, a^{\ast}]$ and $[a^{\ast}, b]$ respectively, there exist $c_1 \in (a, a^{\ast})$ and $c_2 \in  (a^{\ast}, b)$ such that $f^{\sd}(c_2) \leq 0 \leq f^{\sd}(c_1)$. 
  The other remaining cases can be shown in a similar way.
\end{proof}

In order to prove the Quasi-Mean Value Theorem for regular specular derivatives, we compare the graph with lines parallel to its secant line.
Since regular specular derivatives are not linear, the comparison will use the ordinary one-sided derivatives before passing to their specular value.

\begin{theorem} \label{thm:QMVT}
  \emph{(Quasi-Mean Value Theorem)}
  Let $f : [a, b] \to \mathbb{R}$ be a continuous function on $[a,b] \subset \mathbb{R}$.
  Assume $f$ is regularly specularly differentiable in $(a, b)$.
  Then there exist points $c_1$, $c_2 \in (a, b)$ such that 
  \begin{equation*} 
    f^{\sd}(c_2) \leq \frac{f(b)- f(a)}{b-a} \leq f^{\sd}(c_1).
  \end{equation*}
\end{theorem}

\begin{proof}
  Write $A:=\frac{1}{b-a}(f(b)-f(a))$.
  We consider three cases: $f(a) = f(b)$, $f(a) < f(b)$, and $f(a) > f(b)$.
  For starters, suppose $f(a) = f(b)$.
  Define the function $\phi:[a, b]\to\mathbb{R}$ by
  \begin{displaymath}
    \phi(x) = f(x)- f(a).
  \end{displaymath} 
  Clearly, $\phi$ is continuous on $[a,b]$ and regularly specularly differentiable in $(a,b)$.
  Observing that $\phi(a) = \phi(b) = 0$ and $A = 0$, Theorem \ref{thm:Quasi-Rolle} gives points $c_1,c_2\in(a, b)$ such that $\phi^{\sd}(c_2)\le A\le\phi^{\sd}(c_1)$.
  Since $\phi'_\pm = f'_\pm$, Proposition \ref{Prop : Calculating spd} implies $\phi^{\sd} = f^{\sd}$ in $(a, b)$, proving this case.

  Next, assume $f(a) < f(b)$.
  Define the function $\psi :[a,b] \to \mathbb{R}$ by
  \begin{displaymath}
    \psi (x) = A(x-a) + f(a).
  \end{displaymath} 
  If $f=\psi$ on $[a, b]$, then $f'_+ = f'_- = A$ in $(a, b)$, and the conclusion follows.
  Suppose first that there exists $u \in (a,b)$ such that $f(u) > \psi(u)$.
  Choose $k$ with $0 < k < f(u) - \psi(u)$ and define
  \begin{equation*}
    \Psi
    :=
    \{x\in[a, b]:f(x)>\psi(x)+k\}.
  \end{equation*} 
    The set is nonempty, and continuity together with $f(a)=\psi(a)$ and $f(b)=\psi(b)$ gives
  \begin{equation*}
    a < c_1 := \inf\Psi<u<c_2:=\sup\Psi < b
  \end{equation*} 
  and $f(c_i)=\psi(c_i) + k$ for each $i = 1, 2$.
  Thus the two comparison points lie in the interior of the interval.
  There is a sequence $h_j\searrow0$ with $c_1+h_j\in\Psi$.
  For each $j$, we have
  \begin{equation*}
    \frac{f(c_1+h_j)-f(c_1)}{h_j} 
    >
    \frac{\psi(c_1+h_j)+k-\psi(c_1)-k}{h_j}
    =
    A.
  \end{equation*} 
  Since the right-hand derivative exists, $f'_+(c_1)\geq A$.
  For sufficiently small $h > 0$, the definition of $c_1$ also gives
  \begin{equation*}
    \frac{f(c_1)-f(c_1-h)}{h}
    \geq
    \frac{\psi(c_1)+k-\psi(c_1-h)-k}{h}
    =
    A,
  \end{equation*}
  and hence $f'_-(c_1)\geq A$.
  Write $f'_+(c_1)=\tan\theta_1$, $f'_-(c_1) = \tan\theta_2$,  and   $A=\tan\theta_0$, with all three angles in $\left(-\frac{\pi}{2},\frac{\pi}{2}\right)$.
  Then $\theta_1,\theta_2\geq \theta_0$, so Lemma \ref{lem:sd_tan_arctan} implies
  \begin{equation*}
    f^{\sd}(c_1)
    =
    \tan\left(\frac{\theta_1+\theta_2}{2}\right)
    \geq
    \tan\theta_0
    =
    A.
  \end{equation*} 
  At $c_2$, take a sequence approaching from the left within $\Psi$ and use that $f(x)\leq \psi(x)+k$ for $x>c_2$.
  The corresponding difference quotients give $f'_-(c_2)\leq  A$ and $f'_+(c_2)\leq  A$.
  Applying Lemma \ref{lem:sd_tan_arctan}   again yields $f^{\sd}(c_2)\leq  A$.

  If there is no point above $\psi$, there is a point $u$ below it.
  Choose $f(u)-\psi(u)<k<0$ and use $\Psi=\{x\in[a,b]:f(x)<\psi(x)+k\}$.
  The same argument with all inequalities reversed gives the desired bounds after interchanging $c_1$ and $c_2$ .
  Finally, the construction applies without change when $f(a)>f(b)$, since it does not require $A>0$.
 \end{proof}

A regularly specularly differentiable function satisfies the Lipschitz condition whenever its regular specular derivative is bounded.

\begin{corollary}
  Let $f : (a, b) \to \mathbb{R}$ be a continuous function on $(a,b)$.
  Assume $f$ is regularly specularly differentiable and $f^{\sd}$ is bounded on $(a,b)$.
  Let $x_1$, $x_2$ be points in $(a, b)$.
  Then there exists a constant $M > 0$ such that 
  \begin{equation*} 
    \left| f(x_1)-f(x_2) \right| \leq M |x_1-x_2|,
  \end{equation*}
  where $M$ is independent of $x_1$ and $x_2$.
\end{corollary}

\begin{proof} 
  Since $f^{\sd}$ is bounded, $\left|f^{\sd} (x)\right| \leq M$ for some constant $M > 0$ for any $x \in (a,b)$.
  By Theorem \ref{thm:QMVT}, we have
  \begin{equation*} 
    -M \leq \frac{f(x_1)-f(x_2)}{x_1-x_2} \leq M
  \end{equation*}
  for distinct points $x_1,x_2\in(a,b)$, as required.
  The case $x_1=x_2$ is immediate.
\end{proof}

Applying the Quasi-Mean Value Theorem, we show that continuity of $f$ near $x_0$ and the existence of a finite limit of its regular specular derivative at $x_0$ imply classical differentiability at $x_0$.
We make this precise in the following derivative-limit theorem.

\begin{theorem}[Derivative-limit theorem] 
  \label{thm:derivative-limit}
  Let $f:I\to\mathbb{R}$ be continuous, and suppose that $f$ is regularly specularly differentiable on $I\setminus\{x_0\}$, where $x_0\in I$.
  If the limit 
  \begin{equation*}
    \lim_{x\to x_0,\,x\ne x_0}f^{\sd}(x)=\alpha
  \end{equation*}
  exists and is finite, then $f'(x_0)$ exists and $f'(x_0) = f^{\sd}(x_0) = \alpha$.
\end{theorem}

\begin{proof}
  Let $\varepsilon>0$ be given.
  Choose $\delta>0$ such that $(x_0-\delta,x_0+\delta)\subset I$ and
  \begin{equation*}
    \alpha-\varepsilon<f^{\sd}(x)<\alpha+\varepsilon
  \end{equation*}
  whenever $0<|x-x_0|<\delta$.
  Applying Theorem \ref{thm:QMVT} on $[x_0,x_0+h]$ for $0<h<\delta$, there exist $c_1,c_2\in(x_0,x_0+h)$ such that
  \begin{equation*}
    \alpha-\varepsilon<f^{\sd}(c_2) \le\frac{f(x_0+h)-f(x_0)}{h} \le f^{\sd}(c_1)<\alpha+\varepsilon.
  \end{equation*}
  For $-\delta<h<0$, apply the same theorem on $[x_0+h,x_0]$.
  The resulting bounds show that the two-sided difference quotient tends to $\alpha$.
  Hence $f'(x_0)=\alpha$, and Theorem \ref{Thm : ordinary dervatives and specular derivatives} gives $f^{\sd}(x_0)=f'(x_0)$.
\end{proof}

\begin{corollary} \label{cor:cont_of_sd}
  Let $f:I\to\mathbb{R}$ be regularly specularly differentiable on $I$.
  If $f^{\sd}$ is continuous at a point $x \in I$, then $f'(x)$ exists and $f'(x) = f^{\sd} (x)$.
\end{corollary}

\begin{proof}
  Regular specular differentiability implies that $f$ is continuous on $I$.
  Since $f^{\sd}$ is continuous at $x$, we have
  \begin{equation*}
    \lim_{y\to x,\,y\neq x} f^{\sd}(y) = f^{\sd}(x).
  \end{equation*}
  The conclusion follows from Theorem \ref{thm:derivative-limit}.
\end{proof}

Naturally, one can define higher-order regular specular derivatives recursively, as for classical derivatives.
Let $f : I \to \mathbb{R}$  be a function.
Writing $f^{[0]}:=f$ and $f^{[1]}:=f^{\sd}$, for each integer $n \geq 2$, we recursively define the $n$-\emph{th order regular specular derivative} of $f$ at $x \in I$ as
\begin{equation*}
  f^{[n]}(x) := \left(f^{[n-1]}\right)^{\sd}(x)
\end{equation*}
whenever $f^{[n-1]}$ is defined on a neighborhood of $x$ and is regularly specularly differentiable at $x$.
The function $f$ is $n$ times regularly specularly differentiable on $I$ if $f^{[j]}$ is regularly specularly differentiable on $I$ for $j=0,\ldots,n-1$.
In particular, we write the second-order regular specular derivative of $f$ at $x$ as
\begin{displaymath}
  f^{\sd \sd}(x) := \left( f^{\sd} \right)^{\sd} (x).
\end{displaymath} 

\begin{theorem}[Regularity improvement]
  If a function $f:I\to\mathbb{R}$ is twice regularly specularly differentiable on $I$, then $f\in C^1(I)$ and $f'=f^{\sd}$ on $I$.
\end{theorem}

\begin{proof}
  By Definition \ref{def:sd}, $f^{\sd}$ is regularly specularly differentiable on $I$, and hence $f^{\sd}$ is continuous on $I$.
  Corollary \ref{cor:cont_of_sd} gives $f'(x)=f^{\sd}(x)$ for every $x\in I$.
  Thus the classical derivative is continuous on $I$.
\end{proof}

\begin{remark}
  This regularity improvement highlights a key distinction between regular specular differentiability and the broader definition of specular differentiability in \cite{2026a_Jung}.
  Every regularly specularly differentiable function belongs to the class $S^1$ in the terminology of \cite[Prop.~2.20]{2026a_Jung}.
  Two successive regular specular derivatives guarantee continuity of $f^{\sd}$ and hence $C^1$ regularity of $f$.
  However, even a continuous function that is twice specularly differentiable in the sense of \cite{2026a_Jung} need not be classically differentiable; see \cite[Ex.~6.6]{2026a_Jung}.
\end{remark}

\subsection{The Fundamental Theorem of Calculus with specular derivatives}

In this subsection, we establish the Fundamental Theorem of Calculus for regular specular derivatives.
Here, a piecewise continuous function on $[a, b]$ is continuous except at finitely many interior points, at each of which both one-sided limits exist and are finite.
Continuity at the endpoints is understood in the one-sided sense.
For $x\in(a,b)$, we denote the right- and left-hand limits of $f$ by $f(x+)$ and $f(x-)$, respectively.

\begin{theorem}[The Fundamental Theorem of Calculus with specular derivatives, Part 1]
  \label{Thm : FTC with specular derivatives}
  Let $f:[a,b]\to\mathbb{R}$ be piecewise continuous, and define
  \begin{equation*}
    F(x):=\int_a^x f(t)\,dt,
    \qquad x\in[a,b].
  \end{equation*}
  Then $F$ is regularly specularly differentiable on $(a,b)$, and
  \begin{equation}
    \label{eq:FTC-regular-formula}
    F^{\sd}(x)=\mathcal{B}\bigl(f(x+),f(x-)\bigr)
    \qquad\text{for every }x\in(a,b).
  \end{equation}
  At the endpoints, we have $F^{\sd}(a) = F'_+(a)=f(a)$ and $F^{\sd}(b) = F'_-(b) = f(b)$.
\end{theorem}

\begin{proof}
  Since $f$ is bounded and Riemann integrable, $F$ is continuous on $[a,b]$.
  Fix $x\in(a,b)$.
  For sufficiently small $h>0$, we have
  \begin{equation*}
    \frac{F(x+h)-F(x)}{h}
    =\frac{1}{h}\int_x^{x+h}f(t)\,dt.
  \end{equation*}
  Since $f(x+)$ exists, the right-hand side converges to $f(x+)$ as $h\searrow0$.
  Indeed, given $\varepsilon>0$, we have
  \begin{equation*}
    \left|
      \frac{1}{h}\int_x^{x+h}f(t)\,dt-f(x+)
    \right|
    \le
    \frac{1}{h}\int_x^{x+h}|f(t)-f(x+)|\,dt
    \le\varepsilon
  \end{equation*}
  for all sufficiently small $h>0$.
  The value of $f$ at $x$ does not affect the integral.
  Thus $F'_+(x)=f(x+)$.

  Similarly,
  \begin{equation*}
    \lim_{h \searrow 0} \frac{F(x)-F(x-h)}{h}
    =
    \lim_{h \searrow 0} \frac{1}{h}\int_{x-h}^{x}f(t)\,dt
    =
    f(x-),
  \end{equation*}
  and hence $F'_-(x)=f(x-)$. 
  Since both derivatives are finite, $F$ is regularly specularly differentiable at $x$.
  Lemma~\ref{lem:sd_tan_arctan} yields
  \begin{equation*}
    F^{\sd}(x)
    =
    \mathcal{B}\bigl(F'_+(x), F'_-(x)\bigr)
    =
    \mathcal{B}\bigl(f(x+), f(x-)\bigr).
  \end{equation*}
  The same argument at the endpoints, together with the continuity of $f$ there, implies $F'_+(a)=f(a)$ and $F'_-(b)=f(b)$.
\end{proof}

The formula above determines precisely when the regular specular derivative of the integral recovers the integrand at every point.

\begin{corollary}
  \label{cor:pointwise-FTC}
  Let $f$ and $F$ be as in Theorem \ref{Thm : FTC with specular derivatives}.
  The following are equivalent:
  \begin{enumerate}[label=\upshape(\roman*)]
    \item $F^{\sd}(x)=f(x)$ for every $x \in [a, b]$.
    \item $f(s) = \mathcal{B}\bigl(f(s+), f(s-)\bigr)$ for every $s \in (a,b)$ at which $f$ is discontinuous.
  \end{enumerate}
\end{corollary}

\begin{proof}
  The endpoint identities follow immediately from Theorem \ref{Thm : FTC with specular derivatives}.
  At every point $x \in (a,b)$ where $f$ is continuous, Theorem \ref{Thm : FTC with specular derivatives} gives
  \begin{equation*}
    F^{\sd}(x)
    =\mathcal{B}\bigl(f(x+),f(x-)\bigr)
    =\mathcal{B}\bigl(f(x),f(x)\bigr)
    =f(x).
  \end{equation*}
  At every point $x\in(a,b)$ where $f$ is discontinuous, we also have
  \begin{equation*}
    F^{\sd}(x)=\mathcal{B}\bigl(f(x+), f(x-)\bigr).
  \end{equation*}
  Hence $F^{\sd}(x) = f(x)$ holds at every such point if and only if
  \begin{equation*}
    f(x) = \mathcal{B}\bigl(f(x+), f(x-)\bigr)
  \end{equation*}
  at every discontinuity point of $f$.
  This proves the equivalence.
\end{proof}

The following example shows that integration followed by regular specular differentiation does not necessarily recover the original function at its discontinuities.

\begin{example}
  Let $f(x)=\operatorname{sgn}(x)$ on $[-1, 1]$.
  Its indefinite integral is
  \begin{equation*}
    F(x) 
    =
    \int_{-1}^x \operatorname{sgn}(t) \, dt
    =
    |x| - 1.
  \end{equation*}
  At $x=0$, where $f$ has its only jump discontinuity, we have
  \begin{equation*}
    \mathcal{B}\bigl(f(0+), f(0-)\bigr)
    = 
    \mathcal{B}(1, -1) = 0 = f(0).
  \end{equation*}
  Corollary \ref{cor:pointwise-FTC} therefore gives
  \begin{equation*}
    F^{\sd}(x)
    =
    \operatorname{sgn}(x)
    \qquad\text{for every } x \in [-1,1].
  \end{equation*}
  Thus integrating $f$ and then taking the regular specular derivative recovers $f$ at every point of $[-1,1]$.

  Now define
  \begin{equation*}
    g(x)=
    \begin{cases}
      -1 & \text{if } -1 \leq x \leq 0,\\
       1 & \text{if } 0 < x \leq 1.
    \end{cases}
  \end{equation*}
  Since $f$ and $g$ differ only at $0$, their indefinite integrals agree:
  \begin{equation*}
    F(x) 
    = 
    \int_{-1}^x g(t) \, dt 
    = 
    |x| - 1.
  \end{equation*}
  Nevertheless,
  \begin{equation*}
    F^{\sd}(0) = 0 \neq -1 = g(0).
  \end{equation*}
  Thus integrating $g$ and then taking the regular specular derivative does not recover $g(0)$.
\end{example}

We next consider regular specular differentiation followed by integration.

\begin{theorem}[The Fundamental Theorem of Calculus with specular derivatives, Part 2]
  \label{thm:regular-newton-leibniz}
  Let $F:[a,b]\to\mathbb{R}$ be regularly specularly differentiable on $[a,b]$.
  If $F^{\sd}$ is Riemann integrable on $[a,b]$, then
  \begin{equation}
    \label{eq:regular-newton-leibniz}
    \int_a^b F^{\sd}(x)\,dx=F(b)-F(a).
  \end{equation}
\end{theorem}

\begin{proof}
  Write $g:=F^{\sd}$ and let $P=\{a=x_0<x_1<\cdots<x_n=b\}$ be a partition of $[a,b]$.
  Since $g$ is Riemann integrable, it is bounded.
  For each $i=1,\ldots,n$, put
  \begin{equation*}
    m_i:=\inf_{t\in[x_{i-1},x_i]}g(t) 
    \qquad\text{and}\qquad
    M_i:=\sup_{t\in[x_{i-1},x_i]}g(t).
  \end{equation*}
  The function $F$ is continuous on $[a,b]$ by regular specular differentiability.
  Applying Theorem \ref{thm:QMVT} on each interval $[x_{i-1},x_i]$ gives
  \begin{equation*}
    m_i\leq\frac{F(x_i)-F(x_{i-1})}{x_i-x_{i-1}}\leq M_i.
  \end{equation*}
  Multiplying by $x_i-x_{i-1}$ and summing over $i$ yields
  \begin{equation*}
    \sum_{i=1}^{n}m_i(x_i-x_{i-1})
    \leq F(b)-F(a)
    \leq\sum_{i=1}^{n}M_i(x_i-x_{i-1}).
  \end{equation*}
  These are the lower and upper Darboux sums of $g$ for the partition $P$.
  Taking the supremum of the lower sums and the infimum of the upper sums over all partitions, and using the Riemann integrability of $g$, we obtain
  \begin{equation*}
    \int_a^b g(x)\,dx\leq F(b)-F(a)\leq\int_a^b g(x)\,dx.
  \end{equation*}
  This proves \eqref{eq:regular-newton-leibniz}.
\end{proof}

\section{Specular derivatives for multi-variable functions}

Throughout this section, $U$ denotes an open subset of $\mathbb{R}^{n}$, and $\mathbf{e}_i$ denotes the $i$-th standard basis vector of $\mathbb{R}^{n}$.
For $\mathbf{v}=(v_1,\ldots,v_n)\in\mathbb{R}^{n}$ and $s\in\mathbb{R}$, we write $(\mathbf{v},s):=(v_1,\ldots,v_n,s)\in\mathbb{R}^{n+1}$.
We extend regular specular derivatives to functions on $U$ by considering their restrictions to coordinate lines.
The geometric interpretation uses phototangents and tangent hyperplanes in $\mathbb{R}^{n+1}$.

\subsection{Definitions and properties}

Let $f:U\to\mathbb{R}$ be a function, and let $\mathbf{x}\in U$.
For $1\leq i\leq n$, we denote the right- and left-hand partial derivatives of $f$ at $\mathbf{x}$ with respect to $x_i$ by
\begin{equation*}
  \partial^{+}_{x_i}f(\mathbf{x})
  := \lim_{h\searrow0}\frac{f(\mathbf{x}+h\mathbf{e}_i)-f(\mathbf{x})}{h}
  \qquad\text{and}\qquad
  \partial^{-}_{x_i}f(\mathbf{x})
  := \lim_{h\nearrow0}\frac{f(\mathbf{x}+h\mathbf{e}_i)-f(\mathbf{x})}{h},
\end{equation*}
respectively, whenever these limits exist as finite real numbers.

\begin{definition} \label{Def : high-dimensions specularly partial derivatives}
  Let $f:U \to \mathbb{R}$ be a function, and let $\mathbf{x} \in U$.
  We say that $f$ is \emph{regularly specularly partially differentiable} at $\mathbf{x}$ with respect to $x_i$ if both derivatives $\partial^{+}_{x_i}f(\mathbf{x})$ and $\partial^{-}_{x_i}f(\mathbf{x})$ exist.
  We say that $f$ is \emph{regularly specularly differentiable} at $\mathbf{x}$ if $f$ is regularly specularly partially differentiable at $\mathbf{x}$ with respect to $x_i$ for every $i=1,\ldots,n$.
\end{definition}

\begin{example}
  Consider the function $f:\mathbb{R}^{2}\to\mathbb{R}$ defined by
  \begin{equation*}
    f(x,y)=|x+y|.
  \end{equation*}
  Define the set $W:=\{(w_1,w_2)\in\mathbb{R}^{2}:w_1+w_2=0\}$, and let $\mathbf{w}=(w_1,w_2)\in W$.
  Since $f(\mathbf{w})=0$, we have
  \begin{equation*}
    \partial^{+}_{x}f(\mathbf{w})
    =\lim_{h\searrow0}\frac{|w_1+h+w_2|}{h}=1
    \qquad\text{and}\qquad
    \partial^{-}_{x}f(\mathbf{w})
    =\lim_{h\nearrow0}\frac{|w_1+h+w_2|}{h}=-1.
  \end{equation*}
  The same calculation applies to $y$.
  Thus $f$ is regularly specularly partially differentiable with respect to both variables, although neither classical partial derivative exists on $W$.
\end{example}

We define the coordinate phototangents by applying the one-dimensional construction to each coordinate section.
For $\mathbf{x}=(x_1,\ldots,x_n)\in U$ and $1\leq i\leq n$, we denote the coordinate line through $\mathbf{x}$ parallel to the $x_i$-axis by
\begin{equation*}
  \mathbb{R}^{n}_{x_i}(\mathbf{x})
  := \{\mathbf{y}=(y_1,\ldots,y_n)\in\mathbb{R}^{n}:
       y_j=x_j \text{ for all } j\neq i\}.
\end{equation*}

\begin{definition} \label{Def : partial specularly derivatives in high-dimensions}
  Let $f:U\to\mathbb{R}$ be regularly specularly partially differentiable at $\mathbf{x}\in U$ with respect to $x_i$.
  We define the following:
  \begin{enumerate}[label=\upshape(\alph*)]
    \item The \emph{phototangent} of $f$ at $\mathbf{x}$ with respect to $x_i$ is the function
    $\operatorname{pht}_{x_i}f:\mathbb{R}^{n}_{x_i}(\mathbf{x})\to\mathbb{R}$ given by
    \begin{equation*}
      \operatorname{pht}_{x_i}f(\mathbf{y})
      =
      \begin{cases}
        \partial^{+}_{x_i}f(\mathbf{x})(y_i-x_i)+f(\mathbf{x})
        & \text{if } y_i\geq x_i,\\
        \partial^{-}_{x_i}f(\mathbf{x})(y_i-x_i)+f(\mathbf{x})
        & \text{if } y_i<x_i.
      \end{cases}
    \end{equation*}
    \item The (\emph{first-order}) \emph{regular specular partial derivative} of $f$ at $\mathbf{x}$ with respect to $x_i$, denoted by $\partial^{S}_{x_i}f(\mathbf{x})$, is the slope, in the $(x_i,x_{n+1})$-plane, of the line through the two distinct intersection points of the graph of $\operatorname{pht}_{x_i}f$ and the unit sphere $\partial B(\overline{\mathbf{x}},1)$, where $\overline{\mathbf{x}}=(\mathbf{x},f(\mathbf{x}))$.
  \end{enumerate}
\end{definition}

The sphere intersects the two-dimensional plane containing the coordinate phototangent in a unit circle.
Thus Proposition \ref{Prop:sd_iff_pht_cont} ensures that the two intersection points exist and that the slope is finite and independent of the radius.
Equivalently, $\partial^S_{x_i}f(\mathbf{x})=g_i^{\sd}(0)$ for the restriction $g_i(t)=f(\mathbf{x}+t\mathbf{e}_i)$.

\begin{remark}
  Under the assumptions of Definition \ref{Def : partial specularly derivatives in high-dimensions}, Theorem \ref{Thm:sd_criterion} and Lemma \ref{lem:sd_tan_arctan} yield
  \begin{equation} \label{Rmk : specularly partial derivatives formula}
    \partial^S_{x_i}f(\mathbf{x})
    =\lim_{h\searrow0}\mathcal{A}\bigl(g(h),-g(-h),h\bigr)
    =\mathcal{B}\bigl(\partial^+_{x_i}f(\mathbf{x}),\partial^-_{x_i}f(\mathbf{x})\bigr),
  \end{equation}
  where $g(h) = f(\mathbf{x}+h\mathbf{e}_i) - f(\mathbf{x})$.
\end{remark}

We next construct tangent hyperplanes using the intersection points of the coordinate phototangents with the unit sphere.
Recall that a hyperplane in $\mathbb{R}^{n+1}$ is uniquely determined by $n+1$ affinely independent points.

\begin{definition}
  Let $f:U\to\mathbb{R}$ be a function, and let $\mathbf{x}\in U$.
  \begin{enumerate}[label=\upshape(\alph*)]
    \item Let $\mathcal{V}(f,\mathbf{x})$ be the set of indices $i\in\{1,\ldots,n\}$ for which $f$ is regularly specularly partially differentiable at $\mathbf{x}$ with respect to $x_i$.

    \item 
    \label{Def : specularly differentiability for multi-variables (b)}
    Let $\mathcal{P}(f,\mathbf{x})$ be the set of all intersection points of the graphs of $\operatorname{pht}_{x_i}f$, for $i\in\mathcal{V}(f,\mathbf{x})$, with $\partial B(\overline{\mathbf{x}},1)$, where $\overline{\mathbf{x}}=(\mathbf{x},f(\mathbf{x}))$.
  \end{enumerate}
\end{definition}

When no confusion can arise, we write $\mathcal{V}(\mathbf{x})$ and $\mathcal{P}(\mathbf{x})$ for $\mathcal{V}(f,\mathbf{x})$ and $\mathcal{P}(f,\mathbf{x})$.
Writing $\alpha_i=\partial^+_{x_i}f(\mathbf{x})$ and $\beta_i=\partial^-_{x_i}f(\mathbf{x})$, the two points corresponding to $i\in\mathcal{V}(\mathbf{x})$ are
\begin{equation} \label{eq:coordinate-intersections}
  P_i^+=\overline{\mathbf{x}}+\frac{(\mathbf{e}_i,\alpha_i)}{\sqrt{1+\alpha_i^2}}
  \qquad\text{and}\qquad
  P_i^-=\overline{\mathbf{x}}-\frac{(\mathbf{e}_i,\beta_i)}{\sqrt{1+\beta_i^2}}.
\end{equation}
These points are distinct, so $|\mathcal{P}(\mathbf{x})|=2|\mathcal{V}(\mathbf{x})|$.

\begin{example} \label{Ex : not specularly differentiable but specularly partially differentiable}
  Consider the function $f:\mathbb{R}^{2}\to\mathbb{R}$ defined by
  \begin{equation*}
    f(x,y)=
    \begin{cases}
      \displaystyle\frac{x^2}{x^2+y^2} & \text{if } (x,y)\neq(0,0),\\[1em]
      0 & \text{if } (x,y)=(0,0).
    \end{cases}
  \end{equation*}
  At the origin, $\partial^+_y f(0,0)=\partial^-_y f(0,0)=0$, whereas neither one-sided partial derivative with respect to $x$ exists as a finite real number.
  Thus $\mathcal{V}(0,0)=\{2\}$ and $\mathcal{P}(0,0)=\{(0,1,0),(0,-1,0)\}$.
  The function is regularly specularly partially differentiable with respect to $y$, with $\partial^S_y f(0,0)=0$, but is not regularly specularly differentiable at the origin.
\end{example}

\begin{example} \label{Ex : not partially differentiable but specularly partially differentiable}
  Consider the function $f:\mathbb{R}^{2}\to\mathbb{R}$ defined by
  \begin{equation*}
    f(x,y)=
    \begin{cases}
      \displaystyle\frac{xy}{\sqrt{x^2+y^2}} & \text{if } (x,y)\neq(0,0),\\[1em]
      0 & \text{if } (x,y)=(0,0).
    \end{cases}
  \end{equation*}
  The function $f$ is continuous but not classically differentiable at the origin.
  The restrictions of $f$ to both coordinate axes are zero, and hence
  \begin{equation*}
    \partial^+_x f(0,0)=\partial^-_x f(0,0)=\partial^+_y f(0,0)=\partial^-_y f(0,0)=0.
  \end{equation*}
  Thus $f$ is regularly specularly differentiable at $(0,0)$.
  The graphs of the phototangents of $f$ at $(0,0)$ with respect to $x$ and $y$ are the $x$- and $y$-axes in $\mathbb{R}^{3}$, respectively, and
  \begin{equation*}
    \mathcal{P}(0,0)=\{(1,0,0),(-1,0,0),(0,1,0),(0,-1,0)\},
    \qquad
    \partial^S_x f(0,0)=\partial^S_y f(0,0)=0.
  \end{equation*}
\end{example}

\begin{example} \label{ex:coordinate-regularity-not-continuity}
  In one dimension, regular specular differentiability implies continuity, but this implication fails in higher dimensions.
  Consider the function $f:\mathbb{R}^{2}\to\mathbb{R}$ defined by
  \begin{equation*}
    f(x, y)
    =
    \begin{cases}
      \displaystyle\frac{xy}{x^2+y^2} & \text{if } (x,y)\neq(0, 0),\\[1em]
      0 & \text{if } (x,y)=(0,0).
    \end{cases}
  \end{equation*}
  All one-sided partial derivatives at the origin are zero, so $f$ is regularly specularly differentiable at $(0,0)$.
  However, $f(t,t)=\frac{1}{2}$ for every $t\neq0$, so $f$ is not continuous at the origin.
\end{example}

All coordinate phototangents are based at the same point $\overline{\mathbf{x}}=(\mathbf{x},f(\mathbf{x}))$.
To construct a tangent hyperplane that is the graph of an affine function, we must choose $n+1$ intersection points whose projections onto $\mathbb{R}^{n}$ are affinely independent.
The number of intersection points alone does not ensure this condition.

\begin{definition} \label{Def : specular tangent hyperplane}
  Let $f:U\to\mathbb{R}$ be a function, and let $\mathbf{x}\in U$.
  \begin{enumerate}[label=\upshape(\alph*)]
    \item For any $n+1$ points of $\mathcal{P}(f,\mathbf{x})$ whose projections onto $\mathbb{R}^{n}$ are affinely independent, we define the associated \emph{weak specular tangent hyperplane} to be the hyperplane through $\overline{\mathbf{x}}$ parallel to the hyperplane determined by these points.
    We denote the affine function whose graph is this hyperplane by $\operatorname{wstg}f$.

    \item If there is exactly one weak specular tangent hyperplane, we call it the (\emph{strong}) \emph{specular tangent hyperplane} and denote its affine function by $\operatorname{stg}f$.
  \end{enumerate}
\end{definition}

The affine independence of the projected points ensures that each weak specular tangent hyperplane is the graph of an affine function on $\mathbb{R}^{n}$.

\begin{proposition} \label{prop:regular-hyperplane-uniqueness}
  If $f:U\to\mathbb{R}$ is regularly specularly differentiable at $\mathbf{x}\in U$, then the following statements hold.
  \begin{enumerate}[label=\upshape(\alph*)]
    \item 
    \label{prop:regular-hyperplane-uniqueness-1}
    $f$ has at least one weak specular tangent hyperplane.

    \item 
    \label{prop:regular-hyperplane-uniqueness-2}
    $f$ has at most $\binom{2n}{n+1}$ weak specular tangent hyperplanes.

    \item 
    \label{prop:regular-hyperplane-uniqueness-3}
    $f$ has a strong specular tangent hyperplane if and only if all points of $\mathcal{P}(f,\mathbf{x})$ lie in a single hyperplane.
  \end{enumerate}
\end{proposition}

\begin{proof}
  To prove \ref{prop:regular-hyperplane-uniqueness-1}, use the points $P_i^\pm$ in \eqref{eq:coordinate-intersections}.
  Let $\pi:\mathbb{R}^{n+1}\to\mathbb{R}^{n}$ denote the projection onto the first $n$ coordinates.
  For each $j\in\{1,\ldots,n\}$, we have
  \begin{equation*}
    \pi(P_i^+-P_j^-)
    =
    \mathbf{e}_i \left( 1 + \left( \partial^+_{x_i}f(\mathbf{x}) \right)^2 \right)^{-\frac{1}{2}} + \mathbf{e}_j \left( 1 + \left(  \partial^-_{x_j}f(\mathbf{x}) \right)^2 \right)^{-\frac{1}{2}},
    \qquad 1\leq i\leq n.
  \end{equation*}
  Comparing coordinates shows that these $n$ vectors are linearly independent.
  Hence the projections of $P_1^+,\ldots,P_n^+,P_j^-$ are affinely independent.
  Let $H_j$ be the hyperplane determined by these points.
  The parallel translate of $H_j$ through $\overline{\mathbf{x}}$ is a weak specular tangent hyperplane, proving \ref{prop:regular-hyperplane-uniqueness-1}.

  For \ref{prop:regular-hyperplane-uniqueness-2}, recall that $\mathcal{P}(f,\mathbf{x})$ contains $2n$ distinct points.
  Each weak specular tangent hyperplane is obtained from a choice of $n+1$ points with affinely independent projections.
  There are at most $\binom{2n}{n+1}$ such choices, and different choices may determine the same weak specular tangent hyperplane.
  This proves \ref{prop:regular-hyperplane-uniqueness-2}.

  To prove \ref{prop:regular-hyperplane-uniqueness-3}, suppose first that $f$ has a strong specular tangent hyperplane.
  By uniqueness, the parallel translates of all $H_j$ through $\overline{\mathbf{x}}$ coincide, so all $H_j$ are parallel.
  Since every $H_j$ contains $P_1^+$, these hyperplanes coincide.
  Their common hyperplane contains every $P_i^+$ and $P_i^-$, and hence all points of $\mathcal{P}(f,\mathbf{x})$.

  Conversely, suppose that all points of $\mathcal{P}(f,\mathbf{x})$ lie in a single hyperplane $H$.
  Every choice of $n+1$ points with affinely independent projections determines $H$.
  Thus every weak specular tangent hyperplane is the parallel translate of $H$ through $\overline{\mathbf{x}}$.
  Consequently, $f$ has a strong specular tangent hyperplane.
\end{proof}

\begin{remark}
  Suppose that $f$ is regularly specularly differentiable at $\mathbf{x}$ and has a strong specular tangent hyperplane.
  Write the hyperplane containing $\mathcal{P}(f,\mathbf{x})$ as $z=\sum_{i=1}^{n}c_i(y_i-x_i)+d$.
  Since this hyperplane contains $P_i^+$ and $P_i^-$, the slope of their joining line in the $(y_i,z)$-plane is $c_i=\partial^S_{x_i}f(\mathbf{x})$.
  Translating this hyperplane to pass through $(\mathbf{x},f(\mathbf{x}))$ gives
  \begin{equation} \label{Rmk : strong specular tangent hyperplane}
    \operatorname{stg}f(\mathbf{y})
    =
    \sum_{i=1}^{n}\partial^S_{x_i}f(\mathbf{x})(y_i-x_i)+f(\mathbf{x}).
  \end{equation}
  Thus $\operatorname{stg}f(\mathbf{x})=f(\mathbf{x})$ and $\partial_{y_i}(\operatorname{stg}f)(\mathbf{y}) = \partial^S_{x_i}f(\mathbf{x})$ for $\mathbf{y}\in\mathbb{R}^{n}$ and $1\leq i\leq n$, extending Remark \ref{rmk:properties_stl}.
\end{remark}

The following examples show that regular specular differentiability does not imply the existence of a strong specular tangent hyperplane.

\begin{example}
  Consider the function $f:\mathbb{R}^{2}\to\mathbb{R}$ defined by
  \begin{equation*}
    f(x,y)=\left||x|-|y|\right|+|x|+|y|.
  \end{equation*}
  The function $f$ is regularly specularly differentiable at the origin and satisfies $\partial^S_x f(0,0)=\partial^S_y f(0,0)=0$.
  The four points of $\mathcal{P}(f,(0,0))$ lie in the plane $z=\frac{2}{\sqrt{5}}$.
  Hence the strong specular tangent hyperplane of $f$ at the origin is $z=0$.
\end{example}

Here, the following function has nontrivial weak specular tangent hyperplanes.

\begin{example}
  Consider the function $f : \mathbb{R}^2 \to \mathbb{R}$ defined by 
  \begin{equation*}
    f(x, y) = |x|-|y|- x - y.
  \end{equation*}
  The phototangent of $f$ at $(0, 0)$ with respect to $x$ is $\operatorname{pht}_x f : \mathbb{R}^{2}_x(0, 0) \to \mathbb{R}$ by 
  \begin{equation*}
   \operatorname{pht}_x f ( \mathbf{y}_1 ) =
  \begin{cases} 
   -2\mathbf{y}_1 \innerprd \mathbf{e}_1 & \text{if } \mathbf{y}_1 \innerprd \mathbf{e}_1 < 0,\\ 
   0 & \text{if } \mathbf{y}_1 \innerprd \mathbf{e}_1 \geq 0,
  \end{cases}
  \end{equation*}
  for $\mathbf{y}_1 \in \mathbb{R}^2_{x}(0, 0)$ and the phototangent of $f$ at $(0, 0)$ with respect to $y$ is $\operatorname{pht}_y f : \mathbb{R}^{2}_y (0, 0) \to \mathbb{R}$ by
  \begin{equation*}
  \operatorname{pht}_y f ( \mathbf{y}_2 ) =
  \begin{cases} 
   0 & \text{if } \mathbf{y}_2 \innerprd \mathbf{e}_2 < 0,\\ 
   -2\mathbf{y}_2 \innerprd \mathbf{e}_2 & \text{if } \mathbf{y}_2 \innerprd \mathbf{e}_2 \geq 0,
  \end{cases}
  \end{equation*}
  for $\mathbf{y}_2 \in \mathbb{R}^2_{y}(0, 0)$. 
  Note that $f$ is regularly specularly differentiable at $(0, 0)$ with 
  \begin{equation*} 
    \mathcal{P}(0, 0) = \left\{ (1, 0, 0), \left( - \frac{1}{\sqrt{5}}, 0, \frac{2}{\sqrt{5}} \right), (0, -1, 0), \left( 0, \frac{1}{\sqrt{5}}, -\frac{2}{\sqrt{5}} \right) \right\} =: \left\{ \mathbf{p}_1, \mathbf{p}_2, \mathbf{p}_3, \mathbf{p}_4 \right\}.
  \end{equation*}
  For each $i \in \{1,2,3,4\}$, let $\operatorname{wstg}_i f$ denote the affine function whose graph is the weak specular tangent hyperplane associated with the three points $\mathbf{p}_j$, where $j\in\{1,2,3,4\}\setminus\{i\}$.
  Then we see that 
  \begin{align*}
    \displaystyle \operatorname{wstg}_1 f (x, y) &= - \left( \frac{9 - \sqrt{5}}{2} \right)x + \left( \frac{1 - \sqrt{5}}{2} \right)y,\\
    \displaystyle \operatorname{wstg}_2 f (x, y) &= -\left(\frac{1- \sqrt{5}}{2}\right)x + \left(\frac{1- \sqrt{5}}{2}\right)y,\\
    \displaystyle \operatorname{wstg}_3 f (x, y) &= \left( \frac{1 - \sqrt{5}}{2} \right)x - \left( \frac{9 - \sqrt{5}}{2} \right)y,\\
    \displaystyle \operatorname{wstg}_4 f (x, y) &= \left( \frac{1 - \sqrt{5}}{2} \right)x - \left( \frac{1 - \sqrt{5}}{2} \right)y
  \end{align*}
  for $(x, y) \in \mathbb{R}^{2}$.
  Since the four weak specular tangent hyperplanes are distinct, $f$ has no strong specular tangent hyperplane at the origin.
\end{example}

Having $n+1$ intersection points alone does not ensure the existence of a weak specular tangent hyperplane.

\begin{example}
  Consider the function $f:\mathbb{R}^{3}\to\mathbb{R}$ defined by
  \begin{equation*}
    f(x_1,x_2,x_3)=\sqrt{|x_1|}+|x_2|+x_3^2.
  \end{equation*}
  Write $\mathbf{o}=(0,0,0)$.
  Neither one-sided partial derivative with respect to $x_1$ exists as a finite real number, whereas
  \begin{equation*}
    \partial^+_{x_2}f(\mathbf{o})=1,
    \qquad
    \partial^-_{x_2}f(\mathbf{o})=-1,
    \qquad
    \partial^+_{x_3}f(\mathbf{o})
    =
    \partial^-_{x_3}f(\mathbf{o})=0.
  \end{equation*}
  Thus $\mathcal{V}(\mathbf{o})=\{2,3\}$ and
  \begin{equation*}
    \mathcal{P}(\mathbf{o})
    =
    \left\{
      \left(0,\frac{1}{\sqrt{2}},0,\frac{1}{\sqrt{2}}\right),
      \left(0,-\frac{1}{\sqrt{2}},0,\frac{1}{\sqrt{2}}\right),
      (0,0,1,0),
      (0,0,-1,0)
    \right\}.
  \end{equation*}
  Although $|\mathcal{P}(\mathbf{o})| = 4 = n+1$, the four intersection points determine the hyperplane $x_1=0$ in $\mathbb{R}^{4}$.
  Their projections onto $\mathbb{R}^{3}$ lie in the plane $x_1=0$ and are therefore affinely dependent.
  Hence $f$ has no weak specular tangent hyperplane at $\mathbf{o}$ in the sense of Definition \ref{Def : specular tangent hyperplane}.
\end{example}

\subsection{The specular gradient and specular directional derivatives}

We define the regular specular gradient by collecting the regular specular partial derivatives.

\begin{definition}
  Let $f:U\to\mathbb{R}$ be regularly specularly differentiable at $\mathbf{x}\in U$.
  The \emph{regular specular gradient} of $f$ at $\mathbf{x}$ is the vector
  \begin{equation*}
    D^S_{\mathbf{x}}f(\mathbf{x})
    :=
    \left( \partial^S_{x_1}f(\mathbf{x}), \partial^S_{x_2}f(\mathbf{x}), \ldots, \partial^S_{x_n}f(\mathbf{x}) \right).
  \end{equation*}
\end{definition}

When there is no danger of confusion, we write $D^S$ for $D^S_{\mathbf{x}}$.

\begin{remark}
  If $f$ is regularly specularly differentiable at $\mathbf{x}$ and has a strong specular tangent hyperplane at $\mathbf{x}$, formula \eqref{Rmk : strong specular tangent hyperplane} becomes
  \begin{equation*}
    \operatorname{stg}f(\mathbf{y})=D^S f(\mathbf{x})\innerprd(\mathbf{y}-\mathbf{x})+f(\mathbf{x}).
  \end{equation*}
  For a regularly specularly differentiable function $f$, the vector $D^S f(\mathbf{x})$ is defined even when $f$ has more than one weak specular tangent hyperplane at $\mathbf{x}$.
\end{remark}

We consider a simple example of the regular specular gradient.

\begin{example}
  Consider the function $f : \mathbb{R}^{n} \to \mathbb{R}$ defined by 
  \begin{equation*} 
    f( \mathbf{x} ) = \left\vert x_1 \right\vert + \left\vert x_2 \right\vert + \cdots + \left\vert x_n \right\vert, \qquad \mathbf{x} = ( x_1, x_2, \ldots, x_n ).
  \end{equation*}
  For each $i = 1, 2, \ldots, n$, we have $\partial^S_{x_i} f ( \mathbf{x} ) = \operatorname{sgn} ( x_i )$, and hence
  \begin{equation*} 
    D^S f ( \mathbf{x} ) = \left( \operatorname{sgn} ( x_1 ), \operatorname{sgn} ( x_2 ), \ldots, \operatorname{sgn} ( x_n ) \right).
  \end{equation*}
\end{example}

Continuity of the regular specular gradient implies classical differentiability, extending the one-dimensional conclusion of Corollary \ref{cor:cont_of_sd}.

\begin{theorem} \label{thm:continuous-regular-specular-gradient}
  Let $f:U\to\mathbb{R}$ be regularly specularly differentiable on $U$.
  If $D^S f$ is continuous at $\mathbf{x}\in U$, then $f$ is classically differentiable at $\mathbf{x}$ and $Df(\mathbf{x})=D^S f(\mathbf{x})$.
  In particular, if $D^S f$ is continuous on $U$, then $f\in C^1(U)$.
\end{theorem}

\begin{proof}
  Let $\varepsilon>0$.
  Choose $\delta>0$ such that $B(\mathbf{x},\delta)\subset U$ and
  \begin{equation*}
    \left\vert \partial^S_{x_i}f(\mathbf{y}) - \partial^S_{x_i}f(\mathbf{x}) \right\vert  < \varepsilon
  \end{equation*}
  for every $\mathbf{y}\in B(\mathbf{x},\delta)$ and each $i=1,\ldots,n$.
  Let $\mathbf{h}=(h_1,\ldots,h_n)\in\mathbb{R}^{n}$ satisfy $\|\mathbf{h}\|_{\mathbb{R}^n} < \delta$.
  Set $\mathbf{x}^{(0)}=\mathbf{x}$ and define
  \begin{equation*}
    \mathbf{x}^{(i)} = \mathbf{x} + \sum_{j=1}^{i} h_j \mathbf{e}_j
  \end{equation*}
  for each $i = 1, 2, \ldots, n$.
  For each $i=1, 2, \ldots,n$ and every $t$ between $0$ and $h_i$, we have
  \begin{equation*}
    \left\|\mathbf{x}^{(i-1)}+t\mathbf{e}_i-\mathbf{x}\right\|_{\mathbb{R}^n}^2
    =
    \sum_{j=1}^{i-1}h_j^2+t^2
    \leq
    \sum_{j=1}^{n}h_j^2
    =
    \|\mathbf{h}\|_{\mathbb{R}^n}^2
    <
    \delta^2.
  \end{equation*}
  Thus each coordinate segment lies in $B(\mathbf{x},\delta)$.
  On each segment the restriction of $f$ is continuous and regularly specularly differentiable.

  Fix $i\in\{1,\ldots,n\}$ and suppose that $h_i\neq0$.
  Define
  \begin{equation*}
    g_i(t):=f\left(\mathbf{x}^{(i-1)}+t\mathbf{e}_i\right)
  \end{equation*}
  for $t$ between $0$ and $h_i$.
  Apply Theorem \ref{thm:QMVT} to $g_i$ on $[0,h_i]$ if $h_i>0$, and on $[h_i,0]$ if $h_i<0$.
  In either case, there exist $t_1,t_2$ strictly between $0$ and $h_i$ such that
  \begin{equation*}
    \partial^S_{x_i}f\left(\mathbf{x}^{(i-1)}+t_2\mathbf{e}_i\right)
    \leq
    \frac{f\left(\mathbf{x}^{(i)}\right)-f\left(\mathbf{x}^{(i-1)}\right)}{h_i}
    \leq
    \partial^S_{x_i}f\left(\mathbf{x}^{(i-1)}+t_1\mathbf{e}_i\right).
  \end{equation*}
  Since $\mathbf{x}^{(i-1)}+t_k\mathbf{e}_i\in B(\mathbf{x},\delta)$ for each $k=1,2$, the choice of $\delta$ implies
  \begin{equation*}
    \partial^S_{x_i}f(\mathbf{x})-\varepsilon
    <
    \frac{f\left(\mathbf{x}^{(i)}\right)-f\left(\mathbf{x}^{(i-1)}\right)}{h_i}
    <
    \partial^S_{x_i}f(\mathbf{x})+\varepsilon.
  \end{equation*}
  Subtracting $\partial^S_{x_i}f(\mathbf{x})$ and multiplying the resulting absolute-value inequality by $|h_i|$, we obtain
  \begin{equation*}
    \left|
      f\left(\mathbf{x}^{(i)}\right)
      -f\left(\mathbf{x}^{(i-1)}\right)
      -\partial^S_{x_i}f(\mathbf{x})h_i
    \right|
    \leq \varepsilon|h_i|.
  \end{equation*}
  If $h_i=0$, then $\mathbf{x}^{(i)}=\mathbf{x}^{(i-1)}$, so both sides are zero.

  Summing these inequalities and applying the triangle inequality, we obtain
  \begin{equation*}
    \left\vert f(\mathbf{x}+\mathbf{h})-f(\mathbf{x})-D^S f(\mathbf{x})\innerprd\mathbf{h} \right\vert 
    \leq\varepsilon\sum_{i=1}^{n}|h_i|
    \leq\varepsilon\sqrt{n} \|\mathbf{h}\|_{\mathbb{R}^n}.
  \end{equation*}
  Since $\varepsilon$ is arbitrary, $f$ is differentiable at $\mathbf{x}$ with derivative $D^S f(\mathbf{x})$.

  If $D^S f$ is continuous on $U$, this equality holds throughout $U$, proving that $f\in C^1(U)$.
\end{proof}

To define a regular specular directional derivative, we apply the same one-dimensional construction along an arbitrary line.
As Example \ref{ex:coordinate-regularity-not-continuity} shows, the required one-sided directional derivatives need not follow from the existence of the one-sided partial derivatives.

Let $f:U\to\mathbb{R}$ be a function, let $\mathbf{x}\in U$, and let $\mathbf{v}\in\mathbb{R}^{n}$ be a unit vector.
We denote the right- and left-hand directional derivatives of $f$ at $\mathbf{x}$ in the direction $\mathbf{v}$ by
\begin{equation*}
  \partial^+_{\mathbf{v}}f(\mathbf{x})
  :=\lim_{h\searrow0}\frac{f(\mathbf{x}+h\mathbf{v})-f(\mathbf{x})}{h}
  \qquad\text{and}\qquad
  \partial^-_{\mathbf{v}}f(\mathbf{x})
  :=\lim_{h\nearrow0}\frac{f(\mathbf{x}+h\mathbf{v})-f(\mathbf{x})}{h},
\end{equation*}
respectively, whenever these limits exist as finite real numbers.
When all one-sided partial derivatives exist, write
\begin{equation*}
  D^+f(\mathbf{x}) := \left( \partial^+_{x_1}f(\mathbf{x}),\ldots,\partial^+_{x_n}f(\mathbf{x}) \right) 
  \qquad\text{and}\qquad
  D^-f(\mathbf{x}) := \left( \partial^-_{x_1}f(\mathbf{x}),\ldots,\partial^-_{x_n}f(\mathbf{x}) \right).
\end{equation*}

\begin{definition} \label{Def : high-dimensions directional specularly partial derivatives}
  Let $f:U\to\mathbb{R}$ be a function, let $\mathbf{x}\in U$, and let $\mathbf{v}\in\mathbb{R}^{n}$ be a unit vector.
  Suppose that $\partial^+_{\mathbf{v}}f(\mathbf{x})$ and $\partial^-_{\mathbf{v}}f(\mathbf{x})$ both exist.
  For sufficiently small $\delta>0$, define $g:(-\delta,\delta)\to\mathbb{R}$ by
  \begin{equation*}
    g(t):=f(\mathbf{x}+t\mathbf{v}).
  \end{equation*}
  We define the \emph{regular specular directional derivative} of $f$ at $\mathbf{x}$ in the direction $\mathbf{v}$ by
  \begin{equation*}
    \partial^S_{\mathbf{v}}f(\mathbf{x}):=g^{\sd}(0).
  \end{equation*}
\end{definition}

Theorem \ref{Thm:sd_criterion} yields the representation
\begin{equation} \label{def:sdd}
  \partial^S_{\mathbf{v}}f(\mathbf{x})
  =
  \lim_{h \searrow0 } \mathcal{A}\bigl( f(\mathbf{x}+h\mathbf{v})-f(\mathbf{x}), f(\mathbf{x})-f(\mathbf{x}-h\mathbf{v}), h \bigr).
\end{equation}
For $\mathbf{v}=\mathbf{e}_i$, this definition agrees with the regular specular partial derivative: $\partial^S_{\mathbf{e}_i}f=\partial^S_{x_i}f$.
If $f$ is classically differentiable at $\mathbf{x}$, the Chain Rule and Theorem \ref{Thm : ordinary dervatives and specular derivatives} imply
\begin{equation} \label{def:cdd}
  \partial^S_{\mathbf{v}}f(\mathbf{x})
  =
  Df(\mathbf{x})\innerprd\mathbf{v}.
\end{equation}
For a function that is only regularly specularly differentiable, the coordinate partial derivatives need not determine its directional derivatives.

The regular specular directional derivative can instead be calculated directly from its two one-sided directional derivatives.

\begin{proposition} \label{prop:calculation_of_sdd}
  Let $f:U\to\mathbb{R}$ be a function, let $\mathbf{x}\in U$, and let $\mathbf{v}\in\mathbb{R}^{n}$ be a unit vector.
  Suppose that $\alpha:=\partial^+_{\mathbf{v}}f(\mathbf{x})$ and $\beta:=\partial^-_{\mathbf{v}}f(\mathbf{x})$ both exist as finite real numbers.
  Then the regular specular directional derivative $\partial^S_{\mathbf{v}}f(\mathbf{x})$ exists with 
  \begin{equation*}
    \partial^S_{\mathbf{v}}f(\mathbf{x})
    =
    \mathcal{B}(\alpha,\beta)
    =
    \begin{cases}
      \displaystyle\frac{\alpha\beta-1+\sqrt{(\alpha^2+1)(\beta^2+1)}}{\alpha+\beta}
      & \text{if } \alpha+\beta\neq0,\\[1em]
      0 & \text{if } \alpha+\beta=0.
    \end{cases}
  \end{equation*}
  Also, the following statements hold.
  \begin{enumerate}[label=\upshape(\alph*)]
    \item
    \label{prop:calculation_of_sdd-1}
    $\partial^S_{\mathbf{v}}f(\mathbf{x})=0$ if and only if $\alpha+\beta=0$.

    \item
    \label{prop:calculation_of_sdd-2}
    $\operatorname{sgn}(\partial^S_{\mathbf{v}}f(\mathbf{x}))=\operatorname{sgn}(\alpha+\beta)$.

    \item
    \label{prop:calculation_of_sdd-3}
    $\displaystyle \left\vert \partial^S_{\mathbf{v}} f(\mathbf{x}) \right\vert \leq \frac{|\alpha + \beta|}{2}$.

    \item
    \label{prop:calculation_of_sdd-4}
    If $f$ is Lipschitz with constant $L$ on a neighborhood of $\mathbf{x}$, then $|\partial^S_{\mathbf{v}}f(\mathbf{x})|\leq L$ for every direction in which the regular specular directional derivative exists.
  \end{enumerate}
\end{proposition}

\begin{proof}
  Choose $\delta>0$ such that $B(\mathbf{x},\delta)\subset U$ and define $g:(-\delta,\delta)\to\mathbb{R}$ by
  \begin{equation*}
    g(t):=f(\mathbf{x}+t\mathbf{v}).
  \end{equation*}
  Since $g(0)=f(\mathbf{x})$, we have
  \begin{equation*}
    g'_+(0)
    =
    \lim_{h\searrow0}\frac{f(\mathbf{x}+h\mathbf{v})-f(\mathbf{x})}{h}
    =
    \alpha 
    \qquad\text{and}\qquad
    g'_-(0)
    =
    \lim_{h\nearrow0}\frac{f(\mathbf{x}+h\mathbf{v})-f(\mathbf{x})}{h}
    =
    \beta.
  \end{equation*}
  Both derivatives are finite, and hence $g$ is regularly specularly differentiable at $t=0$.
  By the definition of the regular specular directional derivative, $\partial^S_{\mathbf{v}}f(\mathbf{x}) = g^{\sd}(0)$.

  First, apply Lemma \ref{lem:sd_tan_arctan} to $g$ at $t = 0$:
  \begin{equation*}
    \partial^S_{\mathbf{v}}f(\mathbf{x})
    =
    \tan\left(\frac{\arctan\alpha+\arctan\beta}{2}\right)
    =
    \mathcal{B}(\alpha,\beta).
  \end{equation*}
  Applying Proposition \ref{Prop : Calculating spd} with $g'_+(0)=\alpha$ and $g'_-(0)=\beta$ gives the stated piecewise formula.

  The assertions \ref{prop:calculation_of_sdd-1}, \ref{prop:calculation_of_sdd-2}, and \ref{prop:calculation_of_sdd-3} follow from Corollary \ref{cor:calculating_sd} applied to $g$ at $t = 0$.
  To prove \ref{prop:calculation_of_sdd-4}, suppose that $f$ is Lipschitz with constant $L$ on a neighborhood of $\mathbf{x}$.
  Choose $0<r<\delta$ such that the Lipschitz condition holds on $B(\mathbf{x},r)$.
  For all $s,t\in(-r,r)$, we have
  \begin{equation*}
    |g(t)-g(s)|
    =
    |f(\mathbf{x}+t\mathbf{v})-f(\mathbf{x}+s\mathbf{v})|
    \leq
    L\|(t-s)\mathbf{v}\|_{\mathbb{R}^n}
    =
    L|t-s|,
  \end{equation*}
  since $\|\mathbf{v}\|_{\mathbb{R}^n} = 1$.
  Thus $g$ is Lipschitz with constant $L$ on a neighborhood of $0$.
  Applying Corollary \ref{cor:calculating_sd} to $g$ at $t=0$, we obtain $\left|\partial^S_{\mathbf{v}}f(\mathbf{x})\right| = |g^{\sd}(0)| \leq L$.
\end{proof}

\begin{corollary} \label{cor:regular-specular-gradient-estimate}
  Let $f:U\to\mathbb{R}$ be regularly specularly differentiable at $\mathbf{x}\in U$.
  Then
  \begin{equation*}
    \|D^S f(\mathbf{x})\|_{\mathbb{R}^n}
    \leq
    \frac{1}{2}\|D^+f(\mathbf{x})+D^-f(\mathbf{x})\|_{\mathbb{R}^n}.
  \end{equation*}
\end{corollary}

\begin{proof}
  For each $i=1,\ldots,n$, applying Part \ref{prop:calculation_of_sdd-3} of Proposition \ref{prop:calculation_of_sdd} with $\mathbf{v}=\mathbf{e}_i$ yields
  \begin{equation*}
    |\partial^S_{x_i}f(\mathbf{x})|
    \leq
    \frac{1}{2}|\partial^+_{x_i}f(\mathbf{x})+\partial^-_{x_i}f(\mathbf{x})|.
  \end{equation*}
  Squaring and summing over $i$, we obtain
  \begin{equation*}
    \|D^S f(\mathbf{x})\|_{\mathbb{R}^n}^2
    =
    \sum_{i=1}^{n}|\partial^S_{x_i}f(\mathbf{x})|^2
    \leq
    \frac{1}{4}\sum_{i=1}^{n}|\partial^+_{x_i}f(\mathbf{x})+\partial^-_{x_i}f(\mathbf{x})|^2
    =
    \frac{1}{4}\|D^+f(\mathbf{x})+D^-f(\mathbf{x})\|_{\mathbb{R}^n}^2.
  \end{equation*}
  Taking square roots proves the inequality.
\end{proof}

Even when every regular specular directional derivative exists, it need not equal $D^S f(\mathbf{x})\innerprd\mathbf{v}$.

\begin{example} \label{ex:specular-gradient-directional-failure}
  Define the function $f:\mathbb{R}^{2}\to\mathbb{R}$ by
  \begin{equation*}
    f(x,y)=
    \begin{cases}
      \displaystyle\frac{xy(x+y)}{x^2+y^2}
      & \text{if } (x,y)\neq(0,0),\\[1em]
      0 & \text{if } (x,y)=(0,0).
    \end{cases}
  \end{equation*}
  The function $f$ is continuous at the origin.
  Since $f(h,0)=f(0,h)=0$ for every $h\in\mathbb{R}$, we have
  \begin{equation*}
    \partial^+_x f(0,0)=\partial^-_x f(0,0)
    =\partial^+_y f(0,0)=\partial^-_y f(0,0)=0.
  \end{equation*}
  Thus $f$ is regularly specularly differentiable at the origin and
  \begin{equation*}
    D^+f(0,0)=D^-f(0,0)=D^S f(0,0)=(0,0).
  \end{equation*}
  The estimate in Corollary \ref{cor:regular-specular-gradient-estimate} holds with equality, since both sides are zero.

  For each unit vector $\mathbf{v}=(v_1,v_2)$ and every $h\neq0$, we obtain
  \begin{equation*}
    \frac{f(h\mathbf{v})-f(0,0)}{h}
    =
    \frac{v_1v_2(v_1+v_2)}{v_1^2+v_2^2}
    =
    v_1v_2(v_1+v_2).
  \end{equation*}
  Taking the right- and left-hand limits gives
  \begin{equation*}
    \partial^+_{\mathbf{v}}f(0,0)
    =
    \partial^-_{\mathbf{v}}f(0,0)
    =
    v_1v_2(v_1+v_2).
  \end{equation*}
  Proposition \ref{prop:calculation_of_sdd} therefore implies that
  \begin{equation*}
    \partial^S_{\mathbf{v}}f(0,0)
    =
    \mathcal{B}\bigl(v_1v_2(v_1+v_2),v_1v_2(v_1+v_2)\bigr)
    =
    v_1v_2(v_1+v_2).
  \end{equation*}
  In particular, all regular specular directional derivatives exist at the origin.

  For $\mathbf{v}=\frac{1}{\sqrt{2}}(1,1)$, we have
  \begin{equation*}
    \partial^+_{\mathbf{v}}f(0,0)
    =
    \partial^-_{\mathbf{v}}f(0,0)
    =
    \partial^S_{\mathbf{v}}f(0,0)
    =
    \frac{1}{\sqrt{2}}
    \neq
    0
    =
    D^S f(0,0)\innerprd\mathbf{v}
    =
    D^-f(0,0)\innerprd\mathbf{v}
    =
    D^+f(0,0)\innerprd\mathbf{v}
    .
  \end{equation*}
  Thus neither the one-sided nor the regular specular directional derivatives can, in general, be calculated by these inner products.
\end{example}

\subsection{Specular normal vectors}
\label{subsec:specular-normal-vectors}

For a classically differentiable function $f$, the vector $(Df(\mathbf{x}),-1)$ is perpendicular to the tangent hyperplane at $(\mathbf{x},f(\mathbf{x}))$.
We extend this relation to specular tangent hyperplanes and characterize when the regular specular gradient determines a normal vector.

\begin{definition} \label{def:specular-normal-vector}
  Let $f:U\to\mathbb{R}$ have at least one weak specular tangent hyperplane at $\mathbf{x}\in U$.
  A nonzero vector $\mathbf{u}\in\mathbb{R}^{n+1}$ is a \emph{specular normal vector} to the graph of $f$ at $\overline{\mathbf{x}}=(\mathbf{x},f(\mathbf{x}))$ if it is perpendicular to a weak specular tangent hyperplane at $\overline{\mathbf{x}}$.
\end{definition}

Each weak specular tangent hyperplane is the graph of an affine function and therefore has a unique normal vector whose last component is $-1$.
The results below assume regular specular differentiability and characterize when $(D^Sf(\mathbf{x}),-1)$ is a specular normal vector.

\begin{lemma} \label{lem:coordinate-chord-equations}
  Let $f:U\to\mathbb{R}$ be regularly specularly differentiable at $\mathbf{x}\in U$.
  For $1\leq i\leq n$, write $\alpha_i=\partial^+_{x_i}f(\mathbf{x})$, $\beta_i=\partial^-_{x_i}f(\mathbf{x})$, and
  \begin{equation} \label{eq:coordinate-chord-height}
    d_i:=
    \frac{\alpha_i-\beta_i}{\sqrt{1+\alpha_i^2}+\sqrt{1+\beta_i^2}}.
  \end{equation}
  The line through the points $P_i^+$ and $P_i^-$ in \eqref{eq:coordinate-intersections} is given by
  \begin{equation} \label{eq:coordinate-chord-line}
    \begin{cases}
      y_j=x_j & \text{for } j\neq i,\\
      z=f(\mathbf{x})+\partial^S_{x_i}f(\mathbf{x})(y_i-x_i)+d_i.
    \end{cases}
  \end{equation}
  In particular, this line intersects the line $\mathbf{y}=\mathbf{x}$ at $(\mathbf{x},f(\mathbf{x})+d_i)$.
\end{lemma}

\begin{proof}
  Applying \eqref{eq:regular-weighted-formula} to the $i$-th coordinate restriction yields
  \begin{equation*}
    \partial^S_{x_i}f(\mathbf{x})
    =
    \frac{\alpha_i\sqrt{1+\beta_i^2}+\beta_i\sqrt{1+\alpha_i^2}}{\sqrt{1+\alpha_i^2}+\sqrt{1+\beta_i^2}}.
  \end{equation*}
  Consequently,
  \begin{equation*}
    \frac{\alpha_i-\partial^S_{x_i}f(\mathbf{x})}{\sqrt{1+\alpha_i^2}}
    =
    \frac{\partial^S_{x_i}f(\mathbf{x})-\beta_i}{\sqrt{1+\beta_i^2}}
    =
    d_i.
  \end{equation*}
  By \eqref{eq:coordinate-intersections}, the coordinates of $P_i^+$ and $P_i^-$ satisfy \eqref{eq:coordinate-chord-line}.
  Since $P_i^+\neq P_i^-$, equation \eqref{eq:coordinate-chord-line} describes the line through $P_i^+$ and $P_i^-$.
\end{proof}

\begin{theorem}[Criterion for a strong specular tangent hyperplane]
  \label{thm:strong-specular-hyperplane-criterion}
  Let $f:U\to\mathbb{R}$ be regularly specularly differentiable at $\mathbf{x}\in U$, and let $d_i$ be defined by \eqref{eq:coordinate-chord-height}.
  The following are equivalent:
  \begin{enumerate}[label=\upshape(\roman*)]
    \item \label{item:normal-criterion-strong}
    $f$ has a strong specular tangent hyperplane at $\mathbf{x}$.
    \item \label{item:normal-criterion-heights}
    $d_1=d_2=\cdots=d_n$.
    \item \label{item:normal-criterion-gradient}
    The vector $(D^Sf(\mathbf{x}),-1)$ is a specular normal vector at $\overline{\mathbf{x}}$.
  \end{enumerate}
  When these conditions hold, every specular normal vector at $\overline{\mathbf{x}}$ is a nonzero scalar multiple of $(D^Sf(\mathbf{x}),-1)$.
  In particular, $(D^Sf(\mathbf{x}),-1)$ is the unique specular normal vector whose last component is $-1$.
\end{theorem}

\begin{proof}
  Suppose that \ref{item:normal-criterion-strong} holds.
  By Proposition \ref{prop:regular-hyperplane-uniqueness}, all points of $\mathcal{P}(f,\mathbf{x})$ lie in one hyperplane $H$.
  Since $H$ is the graph of an affine function, write its defining equation as
  \begin{equation*}
    z=\sum_{j=1}^{n}c_j(y_j-x_j)+c_0
  \end{equation*}
  for some $c_0,c_1,\ldots,c_n\in\mathbb{R}$.

  For each $i=1,\ldots,n$, the line through $P_i^+$ and $P_i^-$ lies in $H$.
  By Lemma \ref{lem:coordinate-chord-equations}, this line passes through $\bigl(\mathbf{x},f(\mathbf{x})+d_i\bigr)$.
  Substituting this point into the equation of $H$, we obtain
  \begin{equation*}
    f(\mathbf{x})+d_i
    =\sum_{j=1}^{n}c_j(x_j-x_j)+c_0
    =c_0.
  \end{equation*}
  Thus $d_1=\cdots=d_n=c_0-f(\mathbf{x})$, which proves \ref{item:normal-criterion-heights}.

  Conversely, suppose that \ref{item:normal-criterion-heights} holds, and write $d=d_1$.
  By \eqref{eq:coordinate-chord-line}, every $P_i^+$ and $P_i^-$ lies in the hyperplane
  \begin{equation*}
    z=f(\mathbf{x})+D^Sf(\mathbf{x})\innerprd(\mathbf{y}-\mathbf{x})+d.
  \end{equation*}
  Proposition \ref{prop:regular-hyperplane-uniqueness} now implies \ref{item:normal-criterion-strong}.

  Under \ref{item:normal-criterion-strong}, formula \eqref{Rmk : strong specular tangent hyperplane} implies that every point $(\mathbf{y},z)$ on the strong specular tangent hyperplane satisfies
  \begin{equation*}
    z=D^Sf(\mathbf{x})\innerprd(\mathbf{y}-\mathbf{x})+f(\mathbf{x}).
  \end{equation*}
  Consequently,
  \begin{align*}
    (D^Sf(\mathbf{x}),-1)
    \innerprd
    \bigl(\mathbf{y}-\mathbf{x},z-f(\mathbf{x})\bigr)
    &=
    D^Sf(\mathbf{x})\innerprd(\mathbf{y}-\mathbf{x})
    -(z-f(\mathbf{x}))\\
    &=
    D^Sf(\mathbf{x})\innerprd(\mathbf{y}-\mathbf{x})
    -D^Sf(\mathbf{x})\innerprd(\mathbf{y}-\mathbf{x})\\
    &=0.
  \end{align*}
  Since $(D^Sf(\mathbf{x}),-1)\neq\mathbf{0}$, the vector $(D^Sf(\mathbf{x}),-1)$ is normal to the strong specular tangent hyperplane.
  The strong specular tangent hyperplane is also a weak specular tangent hyperplane, and hence \ref{item:normal-criterion-gradient} holds.

  Finally, suppose that \ref{item:normal-criterion-gradient} holds, and write $\mathbf{N}=(D^Sf(\mathbf{x}),-1)$.
  Choose points $Q_0,\ldots,Q_n\in\mathcal{P}(f,\mathbf{x})$ whose projections onto $\mathbb{R}^{n}$ are affinely independent and whose associated weak specular tangent hyperplane has normal vector $\mathbf{N}$.
  Let $H$ be the hyperplane determined by $Q_0,\ldots,Q_n$.
  Since $H$ is parallel to the associated weak specular tangent hyperplane, the vector $\mathbf{N}$ is also normal to $H$.

  For each $k=0,\ldots,n$, both $Q_k$ and $Q_0$ lie in $H$, and hence $\mathbf{N}\innerprd(Q_k-Q_0) = 0$.
  Define $c:=\mathbf{N}\innerprd(Q_0-\overline{\mathbf{x}})$.
  Then
  \begin{equation*}
    \mathbf{N}\innerprd(Q_k-\overline{\mathbf{x}})  
    =
    \mathbf{N}\innerprd(Q_k-Q_0) + \mathbf{N}\innerprd(Q_0-\overline{\mathbf{x}})
    =
    c
  \end{equation*}
  for every $k=0,\ldots,n$.

  We next show that, for each $i=1,\ldots,n$, at least one of $P_i^+$ and $P_i^-$ belongs to $\{Q_0,\ldots,Q_n\}$.
  Fix $i$ and suppose that neither $P_i^+$ nor $P_i^-$ belongs to $\{Q_0,\ldots,Q_n\}$.
  Every $Q_k$ must then belong to a pair $\{P_j^+,P_j^-\}$ with $j\neq i$.
  By \eqref{eq:coordinate-intersections}, the $i$-th coordinate of the projection of every $Q_k$ is $x_i$.
  Thus all $n+1$ projected points lie in the $(n-1)$-dimensional hyperplane $\{\mathbf{y}\in\mathbb{R}^{n}:y_i=x_i\}$.
  This contradicts the affine independence of the projected points.
  Therefore, at least one point from each pair $\{P_i^+,P_i^-\}$ occurs among $Q_0,\ldots,Q_n$.

  Using \eqref{eq:coordinate-intersections} and Lemma \ref{lem:coordinate-chord-equations}, we calculate
  \begin{equation*}
    \mathbf{N}\innerprd(P_i^+-\overline{\mathbf{x}})
    =
    \frac{\mathbf{N}\innerprd(\mathbf{e}_i,\alpha_i)}{\sqrt{1+\alpha_i^2}}
    =
    \frac{\partial^S_{x_i}f(\mathbf{x})-\alpha_i}{\sqrt{1+\alpha_i^2}}
    =
    -d_i
  \end{equation*}
  and 
  \begin{equation*}
    \mathbf{N}\innerprd(P_i^--\overline{\mathbf{x}})
    =
    -\frac{\mathbf{N}\innerprd(\mathbf{e}_i,\beta_i)}{\sqrt{1+\beta_i^2}}
    =
    \frac{\beta_i-\partial^S_{x_i}f(\mathbf{x})}{\sqrt{1+\beta_i^2}}
    =
    -d_i.
  \end{equation*}
  For each $i$, at least one of $P_i^+$ and $P_i^-$ equals some $Q_k$.
  Hence
  \begin{equation*}
    -d_i
    =
    \mathbf{N}\innerprd(Q_k-\overline{\mathbf{x}})
    =
    c.
  \end{equation*}
  Consequently, $d_1=\cdots=d_n=-c$, which proves \ref{item:normal-criterion-heights}.

  When the equivalent conditions hold, denote the unique weak specular tangent hyperplane by $T$.
  Let $\mathbf{u}=(u_1,\ldots,u_{n+1})$ be any nonzero normal vector to $T$.
  Since $T$ passes through $\overline{\mathbf{x}}=(\mathbf{x},f(\mathbf{x}))$, we have
  \begin{equation*}
    \sum_{i=1}^{n}u_i(y_i-x_i)
    +u_{n+1}(z-f(\mathbf{x}))
    =0
  \end{equation*}
  for every $(\mathbf{y},z)\in T$.
  Substituting
  \begin{equation*}
    z-f(\mathbf{x})
    =\sum_{i=1}^{n}\partial^S_{x_i}f(\mathbf{x})(y_i-x_i)
  \end{equation*}
  from \eqref{Rmk : strong specular tangent hyperplane}, we obtain
  \begin{equation*}
    \sum_{i=1}^{n}
    \left(u_i+u_{n+1}\partial^S_{x_i}f(\mathbf{x})\right)
    (y_i-x_i)
    =0
  \end{equation*}
  for every $\mathbf{y}\in\mathbb{R}^{n}$.
  Taking $\mathbf{y}=\mathbf{x}+\mathbf{e}_i$ gives
  \begin{equation*}
    u_i=-u_{n+1}\partial^S_{x_i}f(\mathbf{x})
  \end{equation*}
  for each $i=1,\ldots,n$.
  Therefore,
  \begin{equation*}
    \mathbf{u}
    =-u_{n+1}(D^Sf(\mathbf{x}),-1)
    =-u_{n+1}\mathbf{N}.
  \end{equation*}
  Since $\mathbf{u}\neq\mathbf{0}$, we must have $u_{n+1}\neq0$.
  Conversely, every nonzero scalar multiple of $\mathbf{N}$ is normal to $T$.
  Thus the normal vectors to $T$ are precisely the nonzero scalar multiples of $\mathbf{N}$.
  If the last component of a normal vector is fixed to be $-1$, then $\mathbf{u}=\mathbf{N}$, proving the final assertions.
\end{proof}

For $n=2$, the criterion $d_1 = d_2$ is equivalent to
\begin{equation*}
  (\alpha_1-\beta_1)\left(\sqrt{1+\alpha_2^2}+\sqrt{1+\beta_2^2}\right) -(\alpha_2-\beta_2)\left(\sqrt{1+\alpha_1^2}+\sqrt{1+\beta_1^2}\right)
  =
  0.
\end{equation*}
If all classical partial derivatives of $f$ exist at $\mathbf{x}$, then $\alpha_i=\beta_i$ and $d_i=0$ for every $i$, so the criterion holds.

The vector $(D^Sf(\mathbf{x}),-1)$ is not a specular normal vector when $f$ has several weak specular tangent hyperplanes.

\begin{example} \label{ex:specular-gradient-not-normal}
  Consider the function $f:\mathbb{R}^{2}\to\mathbb{R}$ defined by
  \begin{equation*}
    f(x,y)=\frac{1}{2}(x+|x|)+\frac{1}{2}y+\frac{3}{2}|y|.
  \end{equation*}
  At the origin, we have $\partial^+_x f(0,0)=1$, $\partial^-_x f(0,0)=0$, $\partial^+_y f(0,0)=2$, and $\partial^-_y f(0,0)=-1$.
  Hence
  \begin{equation*}
    d_1=\frac{1}{\sqrt{2}+1}=\sqrt{2}-1  
    \qquad\text{and}\qquad
    d_2=\frac{3}{\sqrt{5}+\sqrt{2}}=\sqrt{5}-\sqrt{2}.
  \end{equation*}
  These values are distinct.
  By Theorem \ref{thm:strong-specular-hyperplane-criterion}, $f$ has no strong specular tangent hyperplane at the origin, and the vector
  \begin{equation*}
    (D^Sf(0,0),-1)
    =
    \left(\sqrt{2}-1,\sqrt{10}-3,-1\right)
  \end{equation*}
  is not perpendicular to any weak specular tangent hyperplane at the origin.
\end{example}

\bibliographystyle{abbrv}
\bibliography{reference}

\end{document}